\documentclass[11pt]{article}

\usepackage{color}
\usepackage{latexsym}
\usepackage{dsfont}
\usepackage{amssymb}
\usepackage{comment}
\usepackage{graphicx}
\usepackage{amsmath, amsfonts,amssymb,theorem,euscript,array,enumerate,amsfonts,mathrsfs}
\usepackage{hyperref}
\usepackage{appendix}

\usepackage{stmaryrd}

\usepackage{mathtools}
\mathtoolsset{showonlyrefs}

\usepackage{comment}

\usepackage{tikz}
\usetikzlibrary{fit,matrix,chains,positioning,decorations.pathreplacing,arrows}

\usepackage{mathtools}
\mathtoolsset{showonlyrefs}

\usepackage[a4paper,margin=1in,heightrounded]{geometry}

\usepackage{algpseudocode}
\usepackage{algorithm}

\usepackage[normalem]{ulem}

\newcommand{\independent}{\protect\mathpalette{\protect\independenT}{\perp}}
\def\independenT#1#2{\mathrel{\rlap{$#1#2$}\mkern2mu{#1#2}}}

\def \Sum{\displaystyle\sum}
\def \Prod{\displaystyle\prod}

\def \b1{\bf{1}}

\def \N{\mathbb{N}}
\def \R{\mathbb{R}}

\def \E{\mathbb{E}}
\def \F{\mathbb{F}}

\def \P{\mathbb{P}}

\def \mrb{\mathrm{b}}

\def \danger{\boldsymbol{D}}

\def\esssup_#1{\underset{#1}{\mathrm{ess\,sup\, }}}

\def\argmin_#1{\underset{#1}{\mathrm{argmin\, }}}
\def\argmax_#1{\underset{#1}{\mathrm{argmax\, }}}

\def \Bc{{\cal B}}

\def \Ec{{\cal E}}
\def \Fc{{\cal F}}

\def \Lc{{\cal L}}

\def \Oc{{\cal O}}

\def\i{{\bf I}}
\def\ta{{\bf T}}
\def\nta{{\bf NT}}
\def\soc{{\bf S}}

\def\bb{{\bf b}}

\def\bc{{\bf c}}

\def \mfb{\mathfrak{b}}

\def\beqs{\begin{eqnarray*}}
\def\enqs{\end{eqnarray*}}
\def\beq{\begin{eqnarray}}
\def\enq{\end{eqnarray}}

\newtheorem{Theorem}{Theorem}[section]
\newtheorem{Definition}{Definition}[section]
\newtheorem{Proposition}{Proposition}[section]

\newtheorem{Remark}{Remark}[section]

\numberwithin{equation}{section}

\begin{document}


\author{
M\'ed\'eric MOTTE
\footnote{LPSM, Universit\'e de Paris  \sf medericmotte at gmail.com
The author acknowledges support of the DIM MathInnov. 
}
\qquad\quad
Huy\^en PHAM
\footnote{LPSM, Universit\'e de  Paris, and CREST-ENSAE, \sf pham at lpsm.paris
The author acknowledges support of the ANR 18-IDEX-001.  This work was also partially supported by the Chair Finance \& Sustainable Development / the FiME Lab (Institut Europlace de Finance)
}
}

\title{Optimal bidding strategies  for digital advertising}

\maketitle

\abstract{
With the emergence of new online channels and information technology, digital advertising  tends to substitute more and more to traditional advertising by offering the opportunity to companies to 
target  the consumers/users that  are really interested by their products or services. We introduce a novel framework for the study of optimal bidding strategies associated to different types of advertising, namely, commercial advertising for triggering purchases or subscriptions, and social marketing for alerting population  about unhealthy behaviours (anti-drug, vaccination, road-safety campaigns). Our  continuous-time models are based on a common framework encoding users online behaviours via their web-browsing at random times, and the targeted advertising auction mechanism widely used on Internet,  the objective being  to efficiently diffuse advertising information by means of digital channels. Our main results are to provide semi-explicit formulas for the optimal value and bidding policy for each of these problems.  
We show  some sensitivity properties of the solution with respect to model parameters, and analyse how the different sources  of digital information  accessible to users including  the social interactions affect the optimal bid  for advertising auctions. We also study how to efficiently combine targeted adver\-tising and non-targeted advertising mechanisms. Finally, some  classes of examples with fully explicit formulas are derived.  
}

\vspace{5mm}


\vspace{5mm}


\vspace{5mm}


\newpage

\section{Introduction}

Through the emergence of new online channels and information technology, targeted advertising plays a growing role in our society and progressively  replaces traditional forms of advertising like newspapers, billboards, etc. Indeed, companies can minimize wasted advertising costs by targeting directly individuals that are potentially interested by the product the advertiser is promoting. Modern targeted media use historical data on internet (cookies) such as tracking online or mobile web activities of consumers.         

 \vspace{1mm}

Optimal control is a suitable mathematical tool for studying advertising problems, and there is already a large literature on this topic. 
In the classical approach, a  dynamical system for the sales process is modeled and the optimisation is performed over the advertising expenditures process.  
We mention the pioneering works by  \cite{Nerlove:1962aa},  \cite{Vidale:1957aa}, and then important papers by Sethi and his collaborators, see  \cite{Feichtinger:1994aa} for an overview of this research field up to the 90s, 
the more recent  paper in \cite{sethi21}, and other references in \cite{RGSethi}, as well as in  the handbook \cite{hand08}. 
We also mention two other works, one about optimal advertising with delay, studied in \cite{2009Gozzi}, and the other \cite{lonzer11} on a model of optimal advertising with singular control.  
 
 \vspace{1mm}

The past decade has seen a growing academic interest in the economic and operations research community for digital advertising. We mention for instance the paper \cite{levmil10} on the design of online advertising markets,  
\cite{yuan14} for a survey on real-time bidding advertising, \cite{goetal21}  for a study of bidding behaviour by learning, \cite{jinetal18} for a multi-agent reinforcement learning algorithm to bid optimisation,  or \cite{baletal14}, \cite{tilletal20} for an optimal bid control problem in online ad auctions, see also 
\cite{choietal20} for a recent literature review on online display advertising.

 \vspace{1mm}

In this paper, we address the following problem. We consider an Agent {\bf A} (company/association) willing to spread some advertising information $\i$ to Users/Individual, like e.g. (i) the existence of a new product,  a new service, or (ii)  the danger of some behaviour (drug, virus, etc). 
These  informations correspond to the following  two types of advertising models that we shall study:   
\begin{enumerate}
\item {\bf Commercial advertising}, modeling situations where informing an individual triggers a reward for the agent.  
We shall consider two types of rewards: {\it purchase-based reward}, corresponding to a punctual payment from the individual to the agent, and {\it subscription-based reward}, corresponding to  a subscription of the individual to a service proposed by the agent, who then receives a regular fee.
\item {\bf Social marketing}, modeling situations where informing an individual cancels a cost continuously perceived by the agent. 
In contrast with  commercial advertising model, the  objective of the agent  is not to make profit but is rather philanthropic.  The aim is  to  change  people's behaviours and to promote social change by raising awareness  about dangers. Classical social marketing campaigns are anti-drugs, vaccination campaigns, road-safety, or  low-fat diet campaigns. 
From the agent's viewpoint, any individual who is not behaving safely is considered to represent a continuous cost to her. 
\end{enumerate}
The issue for Agent {\bf A} is how to diffuse efficiently the information $\i$ by means of ``mo\-dern" online channels  (digital ad, social networks, etc)? With that aim, we propose a continuous-time model for optimal digital advertising strategies, and an important feature is to consider online behaviours of individuals/users that may interact with each other,  and to derive  how  advertising will affect their information states.   
Compared to the classical models that focus directly on macroscopic  variables like sales process controlled by an advertising process, but often without an explicit modelling of the underlying mechanism, 
our approach starts from  a more ``atomic"   level by describing  the individual's behaviours, their social interactions, by an explicit modelling, hence easier to justify from an 
intrinsic point of view.  In particular, we encode the feature of auctions for targeted advertising in our models, which is a crucial component of online advertising. 
The counterpart, in general, is that such microscopic model is often less tractable than classical macro-scale models. In this work, we aim to provide a detailed model with  
reasonably realistic description while keeping it enough tractable in order to obtain explicit solutions.

\vspace{1mm}

Auctions, in targeted advertising, are used to determine which company will have its ad displayed to a given individual. 
There are marketplaces for digital advertising, called Ad exchange, that enable automated buying and selling of ad space.  Each time  an individual browsers through a publisher content (e.g Google, Yahoo, etc), an ad request is sent out for ad space to be viewed, and  the Ad exchange collects the data and information about the viewer via cookies. Then, an auction process is set up in real time where several  advertisers (companies, influencers) declare their bid for ad display, and the highest bidder wins the {\it ad space} by paying a cost according to the first-price or second-price auction rule. 
The long history of auctions, starting from the groundbreaking works of John Nash (\cite{Nash:1951aa}) and later William S. Vickrey (\cite{Vickrey:1961aa}), and their omnipresence on Internet, illustrate the crucial importance of auction theory, also evidenced by the 2020 Nobel prize in economics, awarded to Milgrom and Wilson, for their contributions to auction theory.

\vspace{1mm}

The output of auctions can be quite challenging to predict even in simple frameworks, and as the overall framework of our models is enough complex, we will not  model endogenously  each bidding company (which would turn our optimal control models into games/equilibrium models) and instead assume that at each targeted advertising auction, the exogenous maximal bid from companies other than our agent is a random variable independent from the past, and identically distributed across auctions. This assumption has the practical advantage to keep the control problem tractable.
Finally, one of our models will also encode social interactions allowing individuals who saw the ad to become themselves vectors of information. Again, our modeling of social interactions will be quite simple and symmetric, to keep the problem tractable. For a detailed overview of information spreading models in populations, we refer to \cite{Acemoglu:2011aa}.

 \vspace{1mm}

Besides the different nature of their applications, the both aforementioned studies also differ in their goals. On one hand, commercial targeted advertising is already widely spread on the Internet, and in this case, our study 
proposes a model that could potentially improve  company's bidding strategies. On the other hand, social marketing does not seem currently to use targeted advertising a lot, instead relying more on classical non-targeted advertising, and in this case, our model proposes a method to combine non-targeted advertising with targeted advertising for any organisation or association with that philanthropic purpose.

\vspace{2mm}

\noindent{\bf Our main contributions.} Our first  contribution  is to propose four advertising models, based on a common core framework explicitly modelling on the one hand  individuals  online behaviours via their web-browsing at 
Poisson random times, and on the other hand advertising auctions, each designed for various types of advertising, as described above.  
For each of these problems, we obtain a semi-explicit form of the optimal value function and optimal bidding policy. 


Our second contribution is to propose in one of these models a rich population model, involving individuals spontaneously finding an information,  combining targeted advertising and non-targeted advertising auctions, and highlighting  the role of social interactions. 


Our third contribution is to provide classes of examples where the solutions (optimal value and bidding policy) are fully explicit.


By analysing the form of the solutions, we are able to clearly understand some interesting points:  (1) we observe that the optimal bid to make in a given targeted advertising auction depends not only upon the distribution of other bidders' maximal bids, but also upon the online behaviour of the individual (intensities  at which he connects at random times  to various types of websites);  (2) in the fourth model, involving a population, and adding non-targeted advertising and social interactions in the population, we are able to understand (i) how the presence of social interactions impact the optimal bid to make, and (ii) how the optimal bid to make for non-targeted advertising auctions relates to the optimal bid for targeted advertising auctions and to  the proportion of already informed people. More generally, this work shows 
how the different sources of information (targeted/non-targeted advertising, websites containing the information, and social interactions) affect each other, and in particular how they affect the optimal bids to make in advertising auctions. This is our fourth contribution. 

\vspace{1mm}

The mathematical method for solving  these problems is based on martingale tools, in particular, on techniques involving Poisson processes and their compensators. By means of these techniques and with a suitable change of variable in order to reformulate the problem in terms of proportion of informed individuals, we essentially prove the results in two steps: 1) bounding from above (resp. from below) the optimal value when it is a gain (resp. a cost), and then 2) providing a well chosen policy such that the inequalities in 1) become equalities, thus simultaneously proving that the optimal value is equal to its bound, and obtaining an optimal policy reaching it.

\vspace{1mm}

\noindent{\bf Outline of the paper.} We introduce in Section \ref{sect-core} the core framework of our different  models. In Section \ref{sect-com}, we study two targeted advertising models designed for applications to commercial advertising, the first one modeling advertising to trigger a purchase, the second one modeling advertising to trigger a subscription. In Section \ref{sect-soc}, we study two advertising models applied to social marketing,  
the first one with an arbitrary discount factor, the second one with no discounting, but with extra features of non-targeted advertising and social interactions. We also derive some  insightful properties of the solution related to the 
sensitivity of the optimal bidding strategies with respect to  the individual online behaviour and social interactions effect.   Section \ref{sec:example} presents some  classes of examples where fully explicit formulas can be derived, 
The proofs of our main results are postponed in Section \ref{sec:proofs}, and we conclude in Section \ref{sec:conclusion} by  highlighting some extensions and perspectives.



\section{Basic framework} \label{sect-core}


In this section, we introduce the framework on which all the subsequent models are based, and then  enrich  it in various ways.


The core framework essentially consists in modeling (i) the concept of information, (ii) an individual's online behaviour, (iii)  the targeted advertising auction mechanism, (iv) a targeted advertising bidding strategy, and finally in describing  how these four features combine together to determine the information 
dynamic of an individual.

\subsection{The Information and the Agent}

In this work, all our models will be about some {\it Information}. We shall denote it with a capital ``$\i$'' to emphasize that it is a specific piece of information. It could a priori be any information. Let us give a few examples, further discussed in this work. The Information can be:
\begin{itemize}
\item the existence of a new company,
    \item the existence of a new service (e.g. in Netflix, Amazon, etc),
    \item the existence of a new product (smartphone, computer),
    \item the unhealthiness or healthiness of a behaviour (drug/alcohol consuming, road safety, sexual safety, etc). 
\end{itemize}
In the various models studied in this paper, each model will naturally correspond to one of these three types of information, but for now, let us simply consider a generic Information.

The main characteristic of the Information is that any individual can either {\it not know it} or {\it know it}. In other words, the Information is naturally associated to a binary state for any individual: an individual in state $0$ means that he is not aware of the Information, while an individual in state $1$ means that he is aware of 
the Information.

\vspace{1mm}

In our work, the Agent {\bf A} will represent any entity (company, association, etc) desiring to spread the Information to individuals or population.
\begin{itemize}
    \item In the case of a new service or product, she corresponds naturally to  the company that proposes this service or sells  this product.
    \item In the case of the unhealthiness or healthiness of a given behaviour, she represents  a philanthropic association or a governmental entity aiming to work for social welfare. 
\end{itemize}
The main characteristics of the Agent is that 1) she wants to spread the Information, 2) she has a gain or cost function depending upon how the information spreads, and 3) she will use a digital  advertising strategy as a channel  to diffuse the Information.


\subsection{The Individual and the Action}

Let us start by modelling the general behaviour of an individual. Our model is in continuous time. An individual is associated to some random times when he browses on Internet with the following possible choices:
\begin{itemize}
\item Spontaneously connect to a website providing the Information. Websites intrinsically providing the information are numerous, depending upon the kind of information: specialized websites relaying the Information, company/association's own website, etc. Essentially, any website such that the Information is in the actual website's content, as opposed to the alternative option:
\item Visit a website not  providing {\it a priori}  the information, but displaying targeted ads, and thus susceptible to display the Information  whenever  the agent (company, association, etc) wins the ad auction and pays for it. Important websites displaying targeted ads typically are social networks and search engines.
\end{itemize}


An Individual is associated to independent Poisson processes $(N^\i, N^\ta)$ with respective intensities $\eta^\i$, $\eta^\ta$. $N^\i$ counts the times when the Individual connects to websites intrinsically providing the Information, while  $N^\ta$ counts the times when the Individual connects to websites displaying targeted ads.  
We shall, in our fourth model, introduce a population with several individuals modeled on this basis, each with their own Poisson processes, independent across individuals.

\vspace{2mm}

The Agent aims  to spread  the Information  in order to trigger an {\it Action} from indi\-viduals. The {\it Action} depends upon the type of the Information:
\begin{itemize}
    \item If the Information is about the existence of a service, the expected Action is a {\it subscription}.
    \item If the Information is about the existence of a product, the expected Action is a {\it purchase}.
    \item If the Information is about an unhealthy behaviour, the expected Action is a {\it healthier behaviour}.
\end{itemize}
In this work, we assume that the Agent knows the individuals well enough to be aware of who would do the Action if they had the Information (who would subscribe to the service if he learns that it exists, buy the product if he learns that it exists, or stop some behaviour if he learns that it is unhealthy).

The individuals who would not perform the Action, even informed, are dismissed: the Agent does not try to send them an ad. Therefore, we can assume that the individuals considered in this work are all such that
\beqs 
\text{Getting the Information }\Rightarrow\text{ Doing the Action.}
\enqs

\subsection{The targeted advertising auctions and bidding strategies}

When the Individual connects to websites displaying targeted ads, in reality, many influencers are competing to win the right to display their ads to him. The mechanism used by the website to choose which influencer will display her ad is to make them bid for it. Each influencer has the possibility to propose a bid associated to the Individual's characteristics (intensities of his Poisson processes). This ad emplacement allocation mechanism is what we call {\it targeted advertising auctions}. 


Auctions are complex  to study. They involve several bidders, and are thus part of game theory. The current framework is even more complicated since it is dynamic: an auction is opened each time the Individual connects to a website displaying targeted ads. Our goal is to focus  on providing a {\it strategic} tool to the Agent, and keeping the problem tractable is important in this work.


A good compromise to both take targeted advertising auctions into account while having a strategically solvable problem is to model the maximal bid make by the other bidders (i.e. other than the Agent) as random variables, i.i.d. among auctions. We thus introduce a sequence of i.i.d. real (nonnegative) 
random variables $(B^{\ta}_k)_{k\in\N}$, such that for $k\in\N$, $B^{\ta}_k$ represents the maximal bid of other bidders during the $k$-th targeted advertising auction of the problem.

\vspace{1mm}

We next  introduce the notion of targeted advertising bidding strategies. In essence, a targeted advertising bidding strategy is simply a real valued process $\beta$ which depends {\it at most} from the past, i.e. which cannot depend upon the future (in other words non anticipative), such that at each time $t\in\R_+$, $\beta_t$ represents the bid that the Agent would make if the Individual connects to a website displaying targeted ads.


To rigorously formalize this, let us introduce the filtration $\F$ $=$ $(\Fc_t)_{t\in \R_+}$ generated by the processes $(N^{\i}, N^{\ta}, B^{\ta}_{N^{\ta}})$, i.e.,  
\beqs 
\Fc_t &=& \sigma ((N^{\i}_s, N^{\ta}_s, B^{\ta}_{N^{\ta}_s})_{0\leq s\leq t}), \quad t \geq 0, 
\enqs 
which thus represents all the information about event triggered before time $t$. 


The set of open-loop bidding controls, denoted by $\Pi_{OL}$, is then the set of nonnegative processes $\beta$ predictable and progressively measurable w.r.t. the filtration $\F$.

\subsection{Information dynamic, constant bidding, and advertising cost}
We can now combine all the pieces of modeling previously introduced to define the {\it information dynamic} of the Individual, the notion of constant efficient bidding policy, and the advertising cost. Given an open-loop bidding control $\beta\in\Pi_{OL}$, the information dynamic of the Individual is the $\{0,1\}$-valued process $X^\beta$ satisfying the relation
\begin{equation}
\left\{
\begin{array}{ccl} 
X^\beta_{0^-} &=& 0, \\
dX^\beta_t&=&(1-X^\beta_{t-})(dN^{\i}_t+{\bf 1}_{\beta_t\geq B^{\ta}_{N^{\ta}_t}}dN^{\ta}_t), \quad  t \geq 0. 
\end{array}
\right.
\end{equation}
Let us interpret this dynamic. The individual starts uninformed ($X^\beta_{0^-}=0$). Once he is informed ($X^\beta_t=1$), he stays informed (hence the $(1-X^\beta_{t-})$ part). As long as he is not informed, the remaining part of the dynamic is effective: when the individual connects to a website intrinsically providing the Information, he becomes informed ($dN^\i_t$ part). When he connects to a website displaying targeted ads ($dN^\ta_t$ part), he becomes informed if and only if the Agent's ad is displayed to him, which happens if and only if the Agent wins the auction (${\bf 1}_{\beta_t \geq B^{\ta}_{N^{\ta}_t}}$ part).

\vspace{3mm}

\noindent{\bf Advertising cost.} In the subsequent models, the gain or cost function of the agent will be the combination of 1) a component depending upon the information dynamic of the Individual, and 2) an advertising cost component. The component 1) will depend upon the model, but the advertising cost will always have the same form, namely:  
\begin{align} \label{defCost} 
C(\beta) &= \; \E\Big[\int_0^\infty e^{-\rho t}{\bf 1}_{\beta_t\geq B^{\ta}_{N^{\ta}_t}} \bc(\beta_t,B^{\ta}_{N^{\ta}_t}) d N^{\ta}_t)\Big].
\end{align} 
The interpretation is the following:
\begin{itemize}
\item $\rho\in \R_+$ is a discount rate. Usually, discount rate is chosen to be strictly positive in order to avoid infinite rewards or costs. However, in one of our models (the last one), we will specifically assume $\rho=0$, and it will be an important assumption to make the problem solvable. We shall see that in this model, infinite rewards/costs will never occur despite this assumption.
    \item When the Individual connects to a website displaying targeted advertising ($d N^{\ta}_t$ part), if the targeted advertising auction is won by the agent (${\bf 1}_{\beta_t\geq B^{\ta}_{N^{\ta}_t}}$ part), the agent has to pay a price $\bc(\beta_t,B^{\ta}_{N^{\ta}_t})$, where $\bc:\R^2\rightarrow \R$ is a function depending upon the paying rule defined by the auction. In this paper, the auction payment rule is assumed to be one of the two following standard rules:
    \begin{enumerate}
        \item {\it First-price auctions.} Under this auction rule,  the winner of the auction pays her  bid, and thus, we have $\bc(b,B)=b$.
        \item {\it Second-price auctions.} Under this rule, the winner of the auction pays the {\it second winning bid}, i.e. the bid that she beat. In this case, we have $\bc(b,B)=B$.
    \end{enumerate}
\end{itemize}

\vspace{3mm}

\noindent{\bf Constant bidding policy.} A constant  bidding policy is a constant $b\in\R_+$. The constant bidding control $\beta^b\in\Pi_{OL}$ associated to a constant  bidding policy is defined by the feedback form constraint $\beta^b_t= (1-X^{\beta^b}_{t-})b$. 
It simply models a strategy  where the Agent makes a constant bid $b$ as long as the Individual is not informed (notice that it would be useless to make a positive bid once he is informed).

\vspace{2mm}

We have now introduced all the elements of the core framework. In the sequel, we shall study several advertising problems based on this framework:
\begin{itemize}
    \item In Section \ref{sect-com}, we model commercial advertising problems, i.e. problems where the Agent is a company either trying to sell a service or a product. The common property of both situations is that informing the Individual triggers an Action bringing a {\it reward} to the company (subscription regular fee, purchase punctual fee).
    \item In Section \ref{sect-soc}, we model social marketing problems, i.e. problems where the Agent is an association or government trying to alert people about unhealthy behaviours (anti-drug/alcohol campaigns, road-safety campaigns, etc). The particularity of such type of advertising is that informing people does not bring a reward to the Agent, but instead, it {\it cancels a cost}: as long as an individual has an unhealthy behaviour, he incurs a continuous cost to the philanthropic association. Once informed, he behaves healthier and stops incurring such cost.
\end{itemize}

\section{Commercial advertising model}\label{sect-com}
In this section, we study models for commercial advertising.
The Agent is thus a company trying to maximize its gain.  We will study two types of commercial gains: the subscription-based gain, and the purchase-based gain.

 \subsection{Purchase-based gain function} \label{sec:compur} 
 
We consider the situation where the Information is the existence of a product, where the Agent is a company selling this product, and where the Action of the Individual, once informed, is to purchase the product. We thus  define  the following {\it purchase-based} gain function:
\begin{align} \label{defVpurchase} 
V(\beta) &=\;   \E\Big[\int_0^\infty e^{-\rho t}KdX^\beta_t\Big]-C(\beta), \quad  \mbox{ for } \beta \in\Pi_{OL}. 
\end{align}  
where $K$ is a nonnegative constant. 
Let us interpret this gain function. The part $C(\beta)$ is just the advertising cost from the core framework as defined in \eqref{defCost}. $\rho$ is still the discount rate, and  the part $\int_0^\infty e^{-\rho t}KdX^\beta_t$ simply represents a punctual payment $K$ from the Individual to the Agent when he becomes informed ($dX^\beta_t$ part). This naturally models the reward obtained by the Agent when the individual buys the product. Therefore, $V(\beta)$ represents the net profit of the Agent in the situation of selling a product.

\vspace{2mm}

We now state the main result of this section.

\begin{Theorem}\label{theo-purchase}
We have
\beqs 
V^\star:=\sup_{\beta\in\Pi_{OL}} V(\beta)&=& 
\sup_{b\in\R_+} V(\beta^b),
\enqs 
with
\begin{align} \label{Vbetab} 
V(\beta^b) &= \; \frac{\eta^{\i}K+\eta^{\ta}\E\big[\big(K-\bc(b, B_1^{\ta})\big) {\bf 1}_{b\geq B^{\ta}_1}\big]}{\eta^{\i}+ \rho + \eta^{\ta}\P[b\geq B^{\ta}_1]}, \quad \forall b\in\R_+.
\end{align} 
Furthermore, any $b^\star$ $\in$ 
$\argmax_{b\in\R_+}V(\beta^b)$
yields an optimal constant bid policy, i.e. an optimal open-loop bid  control  $\beta^{b^\star}$  taking the form of a constant  bid.
\end{Theorem}

\paragraph{Interpretation.} Let us interpret this result  by first understanding the role of $\rho$. It is well known that a discount rate is mathematically equivalent to a random termination date of the problem following an exponential distribution $\Ec(\rho)$  with parameter $\rho$. 
Up to adding this random termination time, we can thus consider that the problem has no discount rate. Given this interpretation, and assuming that the Agent plays a constant bidding policy $b$, notice that the inner fraction in $V(\beta^b)$  can be written as
\beqs 
\pi_{\i}K+\pi_{\rho}\times 0 + \pi_{\ta}\E [K-\bc(b, B_1^{\ta})\mid b\geq B^{\ta}_1]
\enqs 
where $(\pi_{\i}, \pi_{\rho},\pi_{\ta})$ are probability weights proportional to $(\eta^{\i},\rho , \eta^{\ta}\P[b\geq B^{\ta}_1])$. This expression should be seen as the expected reward of the Agent computed in terms of how the problem terminates:
\begin{itemize}
    \item When it  terminates with the Individual finding the Information by himself with probability $\pi_{\i}$, the Agent only perceives the reward $K$.
    \item When it  terminates at the random time associated to $\Ec(\rho)$,  the Individual has not had the time to be informed: the Agent perceives nothing.
    \item When it  terminates with the Individual getting informed by viewing the Agent's targeted ad with probability $\pi_{\ta}$,  the Agent perceives $K$ and pays $\bc(b, B^{\ta}_1)$ because he had to pay the auction's price.
\end{itemize}


\vspace{1mm}

Besides the quantitative aspect of this result, an important qualitative property is that a constant bidding policy is enough to reach the optimal value over all open-loop bidding controls. This is particularly interesting from a model-free viewpoint (reinforcement learning) as it means that one can restrict the search for an optimal strategy to the set of nonnegative constant bidding policies, 
which is a reasonably ``small'' set.


\begin{Remark}[Cost dual viewpoint] 
Another interesting way to formulate the optimal value and bid is from a {\it cost viewpoint} (and this is actually how we prove this formula in this paper): the idea is to consider the best possible scenario for the Agent, which arises when  the Individual directly connects to a website containing the information from the very beginning, and then look at the real scenario {\it relatively} to this best scenario. The real scenario necessarily brings a smaller gain than the best scenario, and thus, it is {\it as if} the Agent won the best scenario gain but then pays  a cost corresponding to this difference. From this viewpoint, the goal is to minimize this cost. The best scenario gain  is clearly equal to  $K$, and we can rewrite  the optimal value from \eqref{Vbetab} as 
\beqs 
\sup_{\beta\in\Pi_{OL}} V(\beta)&=& K- 
\inf_{b\in \R_+} 
\frac{\rho K+\eta^{\ta}\E\big[\bc(b, B_1^{\ta}) {\bf 1}_{b\geq B^{\ta}_1}\big]}{\eta^{\i}+ \rho + \eta^{\ta}\P[b\geq B^{\ta}_1]},
\enqs 
and any $b^\star\in\R_+$ such that 
\beqs 
b^\star&=&
\argmin_{b\in\R_+} 
\frac{\rho K+\eta^{\ta}\E\big[\bc(b, B_1^{\ta}) {\bf 1}_{b\geq B^{\ta}_1}\big]}{\eta^{\i}+ \rho + \eta^{\ta}\P[b\geq B^{\ta}_1]}
\enqs 
yields an optimal constant bid.
\end{Remark}

\vspace{1mm}

We next introduce the following special optimal minimal bidding policy.

\begin{Definition}[Smallest optimal constant bid policy]
We denote $b^\star_{min}$ the constant bidding policy such that
\beqs 
b_{min}^\star &=& \min \argmax_{b\in\R_+}V(\beta^b).
\enqs 
$b_{min}^\star$ is called the smallest optimal constant bidding policy.
\end{Definition}

\begin{Remark}
From the proofs of our results, it is possible to see that the open-loop bidding control $\beta^{b_{min}^\star}$ is the smallest optimal open-loop bidding control, i.e. for all optimal open-loop bidding control $\beta$, we have 
$\beta^{b_{min}^\star}_t$ $\leq$ $\beta_t$, $(\omega,t)$-a.e. 
\end{Remark}

We have the following result about the sensitivity to parameters and upper bounds of the optimal value and smallest minimal optimal bidding policy.

\begin{Proposition}\label{prop-sensitivity}
The optimal value $V^\star$ is increasing in $\eta^{\i},\eta^{\ta}$ and decreasing in $\rho$, and the smallest optimal constant bidding policy $b_{min}^\star$ is decreasing in $\eta^{\i},\eta^{\ta}$ and increasing in $\rho$. 
Finally, we have
\beqs 
V^\star\geq \frac{\eta^{\i} K}{\eta^{\i}+\rho}, \quad \quad b_{min}^\star \; \leq \;  K-V^\star \; \leq \;  \frac{\rho K}{\eta^{\i}+\rho}.
\enqs 
\end{Proposition}

\vspace{1mm}

The interpretation of the above Proposition is the following. 

\vspace{1mm}

\noindent{\bf Properties for $V^\star$.}  The higher $\eta^{\i}$ is, the more frequently the Individual connects to a website containing the Information, and thus the sooner he would learn by this channel  the Information, which, at fixed constant bid $b$, naturally increases the expected gain of the Agent. This explains why the gain with any constant bid $b$, and thus the optimal gain, of the Agent is increasing in $\eta^{\i}$. Given $\tilde{\eta}^{\ta}\leq \eta^{\ta}$, the Agent can always ``emulate'' any scenario associated to a constant bid $b$ and a frequency $\tilde{\eta}^{\ta}$ of connection to a website displaying targeted ads, simply by constraining himself to bid $b$ only with probability $\frac{\tilde{\eta}^{\ta}}{\eta^{\ta}}$ and $0$ otherwise at each auction: by a standard property of Poisson processes, it will be equivalent to always bid $b$ at auctions occurring with intensity $\tilde{\eta}^{\ta}$. Consequently, with the intensity $\eta^{\ta}$, the Agent can replicate all the gains that an intensity $\tilde{\eta}^{\ta}$ could yield, and thus his optimal gain is increasing in $\eta^{\ta}$. Finally, the larger $\rho$ is, the more impatient the Agent is, and thus the less value he gives to potential future rewards, which explains why his optimal gain $V^\star$ is decreasing in $\rho$. The lower bound of $V^\star$ simply corresponds to the gain associated to the constant bidding policy consisting in bidding $0$ at each auction, i.e. never displaying any ad, and thus simply waiting that the individual informs himself on a website containing the Information.

\vspace{1mm}

\noindent{\bf Properties  for the smallest optimal constant bid.}  When an auction is opened, two scenarii can occur: 1) the Agent wins the auction, receives $K$ and pays the auction price, or 2) the Agent loses the auction, and the problem of informing the individual keeps going, with an optimal value $V^\star$. In other words, if we put the auction price apart, an auction can be seen as providing a reward $K$ when it is won, and $V^\star$ when it is lost. Notice that it is equivalent to consider that $V^\star$ is won anyway, and that the auction is a standard static auction providing the additional reward $K-V^\star$ if the auction is won, and $0$ if it is lost. The larger $V^\star$ is, the smaller $K-V^\star$ is, and thus the smaller the bid that the Agent should be willing to make to win $K-V^\star$ is. The smallest optimal bid is thus decreasing in $V^\star$, which explains why its sensitivity to all the parameters are reversed w.r.t. the sensitivity of $V^\star$. In such auction, it is also clear that the Agent would have no interest in paying more than $K-V^\star$ to win the auction, which justifies the upper bound $K-V^\star$. The greater upper bound $\frac{\rho K}{\eta^{\i}+\rho}$ directly  comes from the lower bound of $V^\star$. In particular, this implies that for getting an optimal bidding constant strategy, we can restrict the search of  the supremum over $b$ in $V(\beta^b)$ to the bounded interval $[0,\frac{\rho K}{\eta^{\i}+\rho}]$.

 \subsection{Subscription-based gain function}
 
 We now consider  the situation where the Information is the existence of a service, where the Agent is the company proposing this service, and where the Action of the Individual, once informed, is to subscribe to the service. To that aim, we then simply  consider the following {\it subscription-based} gain function:
\begin{align} \label{defVsus} 
V(\beta) &= \;  \E\Big[\sum_{n\in\N} e^{-(\tau^\beta+n)\rho}K\Big]-C(\beta), \quad \mbox{ for } \beta \in\Pi_{OL},
\end{align} 
where $\tau^\beta:= \inf\{t\in\R_+: X^\beta_t=1\}$ is the time of information of the individual.

\vspace{3mm}

Let us interpret this gain function.  Again, the part $C(\beta)$ is the advertising cost described in the core framework, and  $\rho$ is still the discount rate. 
The other part 
represents the gain coming from the Individual's information dynamic. It  corresponds to a regular payment of $K$ every period $1$ from the time of information $\tau^\beta$ (and thus the time of subscription) of the Individual.

\vspace{2mm}

 We can now state the main  result of this section. 

\begin{Theorem} \label{theo-sus}
We have
\beqs 
\sup_{\beta\in\Pi_{OL}}  V(\beta) &=&  
\sup_{b\in\R_+} 
V(\beta^b),
\enqs
with
\beqs
V(\beta^b) &=& \frac{\eta^{\i}\frac{K}{1- e^{-\rho}}+\E\big[\big(\frac{K}{1-e^{-\rho}}-\bc(b,B_1^{\ta})\big) {\bf 1}_{b\geq B^{\ta}_1}]}{\eta^\i+ \rho + \P[b\geq B^{\ta}_1]}, 
\enqs 
and any $b^\star$ $\in$ 
$\argmax_{b\in\R_+} V(\beta^b)$ 
yields an optimal constant bid, i.e. an optimal open-loop bid control taking the form of a constant bid.
\end{Theorem}

\noindent{\bf Interpretation.} Notice that the regular payment of $K$ every period of duration $1$ from the time of information is, from the Agent's viewpoint, equivalent to a unique payment of 
$\frac{K}{1-e^{-\rho}}$ at the time of information. We are thus reduced to the previous case of purchase-based gain.

\section{Social marketing models}\label{sect-soc}

We now model a quite  different type of advertising, called {\it social marketing}. Social marketing is the activity of making advertising campaigns not to make profit but to alert people, in particular about unhealthy behaviours (anti-drug campaigns, road-safety campaigns, sexual-safety campaigns, etc). 
The Agent, here, is either a philanthropic  association or a governmental entity working for social welfare, and considers that each Individual not behaving healthily incurs a cost to her. 
As opposed to commercial advertising from previous section, informing an Individual here does not bring a reward to the Agent, but instead, cancels a cost.

\vspace{2mm}



For this application, our study will be split in two sub-cases:
\begin{enumerate}
    \item The case with a positive discount rate $\rho$, based on the same framework as previous models but with a cost function, and
    \item The important case with no discount  rate (i.e. $\rho$ $=$ $0$), where we shall be able to enrich the basic  framework  by introducing a population of $M$ individuals as well as a non-targeted advertising mechanism, therefore turning the model into a population control problem.
\end{enumerate}
In both cases the Agent's goal will be to {\it minimize} her  cost function.

\subsection{Case with a discount rate} \label{secsocialdiscount}

We start by the  case, with no social interaction nor non-targeted advertising, but with an arbitrary discount rate $\rho$. Besides the processes $N^\i$ and $N^\ta$, we consider a third Poisson process 
$N^{\danger}$ ($\danger$ for ``Dangerous behaviour'') , independent from the others, with normalized intensity $\eta^{\danger}=1$, counting the times when the Individual behaves unsafely. 

In this social marketing problem, the cost function of the Agent is defined by
\begin{align} \label{defVsocialdis}  
V(\beta) &= \;  \E\Big[\int_0^\infty e^{-\rho t}K(1-X^\beta_{t-})dN^{\danger}_t\Big]+C(\beta), \quad  \mbox{ for } \beta \in\Pi_{OL}.  
\end{align} 
The part $C(\beta)$ is the advertising cost, and the part $\E\Big[\int_0^\infty  e^{-\rho t}K(1-X^\beta_{t-})dN^{\danger}_t\Big]$  measures the (discounted) cost perceived in the period before the Individual was informed, assuming that the Individual incurs a cost $K$ to the Agent every time he behaves unsafely. 

\vspace{2mm}

We have the following result.

\begin{Theorem} \label{theo-socialdis}
We have
\beqs 
V^\star \; := \;  \inf_{\beta \in\Pi_{OL}}  V(\beta) &=&
\inf_{b\in\R_+}  V(\beta^b),
\enqs
with
\beqs
V(\beta^b) &=& \frac{K+\eta^{\ta}\E\big[\bc(b,B_1^{\ta}) {\bf 1}_{b\geq B^{\ta}_1}\big]}{\eta^{\i}+ \rho + \eta^{\ta}\P[b\geq B^{\ta}_1]}
\enqs 
and any $b^\star$ $\in$ 
$\argmin_{b\in\R_+} V(\beta^b)$  
yields an optimal constant bid, i.e. an optimal open-loop bid taking the form of a constant bid.
\end{Theorem}

\noindent{\bf Interpretation.} Here again, we interpret $\rho$ as the parameter of a random terminal time with exponential distribution. Notice that in the case of social marketing, 
there is already a random terminal time: the time when the Individual connects on the website intrinsically containing the information. Indeed, in such case, the cost stops and the problem stops as well. Both terminal times are exponential random variables with respective parameters $\eta^{\i}$ and $\rho$. It is known that they can be compressed in a unique terminal time (the minimum of both) with parameter $\eta^{\i}+\rho$. In other words, up to replacing the original intensity $\eta^{\i}$ of connection to a website containing the Information by $\eta^{\i}+\rho$, we are reduced to a problem with no discount rate ($\rho=0$). The inner fraction can be split as follows:
\begin{align} \label{decsocial} 
\frac{K+\eta^{\ta}\E\big[\bc(b,B_1^{\ta}) {\bf 1}_{b\geq B^{\ta}_1}\big]}{\eta^{\i}+ \rho + \eta^{\ta}\P[b\geq B^{\ta}_1]} &=\;
\frac{K}{\eta^{\i}+ \rho + \eta^{\ta}\P[b\geq B^{\ta}_1]} 
+\frac{\eta^{\ta}\E\big[\bc(b,B_1^{\ta}) {\bf 1}_{b\geq B^{\ta}_1}\big]}{\eta^{\i}+ \rho + \eta^{\ta}\P[b\geq B^{\ta}_1]},
\end{align} 
and has the following interpretation. $\eta^{\i}+ \rho + \eta^{\ta}\P[b\geq B^{\ta}_1]$ is the intensity of the time of information of the Individual, and thus $\frac{1}{\eta^{\i}+ \rho + \eta^{\ta}\P[b\geq B^{\ta}_1]}$ is the expected time before information. During this time, a continuous cost $K$ is essentially perceived, which explains the  first term in the r.h.s. of \eqref{decsocial}. The second term 
is essentially  the expected cost perceived at the time of termination of the problem, given that in this case, no reward, and only the ad cost, is paid.

\vspace{1mm}

As in the commercial advertising case, we introduce the following special optimal minimal bidding policy.
\begin{Definition}[Smallest optimal constant bid policy]
We denote $b_{min}^\star$ the constant bidding policy such that
\beqs 
b_{min}^\star &=& \min \argmin_{b\in\R_+}V(\beta^b).
\enqs 
$b_{min}^\star$ is called the smallest optimal constant bidding policy.
\end{Definition}

We have the following result about the sensitivity to parameters and upper bounds of the optimal value and smallest minimal optimal bidding policy.

\begin{Proposition}\label{propsocial-sensitivity}
The optimal value $V^\star$ and the smallest optimal bid are decreasing in $\eta^{\i},\eta^{\ta}$, and $\rho$, and we have
\beqs 
V^\star \; \leq \;  \frac{K}{\eta^{\i}+\rho}, \quad \quad b_{min}^\star  \; \leq \;  V^\star \; \leq \;  \frac{K}{\eta^{\i}+\rho}.
\enqs 
\end{Proposition}

\vspace{1mm}

The interpretation of the above proposition is the following.

\vspace{1mm}

\noindent{\bf Properties  for $V^\star$.}  The higher $\eta^{\i}$ is, the more frequently the Individual connects to a website containing the Information, and thus the sooner he would learn by this channel the Information, and stop inducing a cost to the Agent, which naturally decreases the expected cost of the Agent. The justification of the sensitivity in $\eta^{\ta}$ is similar to the corresponding interpretation for commercial advertising. 
Finally, the larger $\rho$ is, the more impatient the Agent is, and thus the less value she gives to potential future costs, which explains why her optimal cost $V^\star$ is decreasing in $\rho$. The lower bound of $V^\star$  corresponds to the cost associated to the constant bidding policy consisting in bidding $0$ at each auction, and thus  waiting that the individual informs himself on a website containing the Information.

\vspace{1mm}

\noindent{\bf Properties  for the smallest optimal constant bid.} When an auction is opened, two scenarii can occur: 1) the Agent wins the auction, and the cost stops, or 2) the Agent loses the auction, and the problem of informing the individual keeps going on, with an optimal cost $V^\star$. In other words, if we put the auction price apart, an auction can be seen as incurring a cost $0$ when it is won, and $V^\star$ when it is lost. Notice that it is equivalent to consider that the cost $V^\star$ is incurred anyway, and that the auction is a standard static auction providing the compensating reward $V^\star$ if the auction is won, and $0$ if it is lost. Thus, the larger $V^\star$ is,  the greater the bid that the Agent should be willing to make to win $V^\star$ is. The smallest optimal bid is thus increasing in $V^\star$, which explains why its sensitivity to all the parameters are the same as $V^\star$. In such auction, it is also clear that the Agent would have no interest in paying more than $V^\star$ to win the auction, which justifies the upper bound $V^\star$. The greater upper bound $\frac{K}{\eta^{\i}+\rho}$  follows  from the upper bound of $V^\star$.

\subsection{Case with no discount rate, with social interactions and non-targeted advertising} \label{sec:nodiscount}

In this section, we consider a social marketing model with no discounting, but with much more features than previous models. 
We do not simply model websites intrinsically containing the Information and websites displaying targeted ads, but also model the alternative for users  to  connect on website displaying non-targeted ads, and to socially interact with each other. 
The arguments for introducing  these two extra features is twofold:
\begin{enumerate}
    \item {\it For relevance in terms of applications.} Social marketing nowadays still widely happens via non-targeted advertising (TV campaigns, etc). Although our model proposes to use targeted advertising, it thus seems important to not completely dismiss the current method, and instead propose a way to combine both mechanisms.
    \item {\it Mathematical reason.} The absence of discount rate  allows the problem to still be tractable even by adding these features. 
\end{enumerate}

Let us reintroduce for sake of completeness and self-contained reading   each component of the framework together  with these additional features.

\paragraph{The population.}

We now  consider a population with $M$ individuals, with online behaviour  characterised by:
\begin{itemize}
    \item a family of $M$ i.i.d. triplets $(N^{m,\i}, N^{m,\ta}, N^{m,\nta}, N^{m,\danger})$,  for $m\in \llbracket 1,M\rrbracket$, where $N^{m,\i}$, $N^{m,\ta}$, $N^{m,\nta}$, and $N^{m,\danger}$ are four independent Poisson processes with respective intensities $\eta^{\i}$, $\eta^{\ta}$, $\eta^{\nta}$, and $\eta^{\danger}=1$.  
    Notice that we assume that the population is homogeneous, i.e.  each individual shares  the same intensities. 
    \item a family $(N^{m,i, \soc})_{m,i\in \llbracket 1,M\rrbracket}$ of i.i.d. Poisson processes with intensity $\eta^{\soc}$, independent from the other Poisson processes.
\end{itemize}
For all $m\in \llbracket 1,M\rrbracket$, the processes $N^{m,\i}$, $N^{m,\ta}$, and $N^{m,\danger}$, have the same interpretation as in the previous model: $N^{m,\i}$ counts the times when individual $m$ visits a website intrinsically containing the Information (in this case, it would be an association's website, the website specialized in health, etc). $N^{m,\ta}$ counts the times when individual $m$ connects to a website displaying targeted ads, and $N^{m,\danger}$ counts the time when he behaves unsafely. The new features are: $N^{m,\nta}$,  counting the times when individual $m$ visits a website displaying {\it non}-targeted ads, and for $m,i\in \llbracket 1,M\rrbracket$, $N^{m,i,\soc}$ counting  the social interactions between individuals $m$ and $i$ in the population.

\paragraph{Targeted and non-targeted advertising auctions.}
\begin{itemize}
\item {\it Targeted advertising auctions.}  For each  individual $m\in \llbracket 1,M\rrbracket$, whenever  he  connects to a website displaying targeted ads, an auction is automatically opened where several agents bid to win the right to display their ads to the individual. As in previous models, 
we model the maximal bid from other bidders (other than our Agent), by introducing 
an i.i.d. family of nonnegative random variables $(B^{m,\ta}_k)_{k\in\N,m\in \llbracket 1,M\rrbracket}$, where $B^{m,\ta}_k$ represents the maximal bid from other bidders at the $k$-th targeted advertising auction concerning individual $m$. We also denote by $\bc^{\ta}:\R^2\rightarrow \R$  the paying rule function of this targeted auction, which is again assumed to be either of first-price or second-price auction rule. 
\item {\it Non-targeted advertising auctions.} In this model, we also consider non-targeted advertising. Each time when an individual (regardless of his index) connects to a website displaying non-targeted ads, here again, agents will compete to display their ads (with the only difference that they cannot make their bid depending upon the individual who connects to the website, hence the name ``non-targeted advertising''). 
An auction is thus also opened at each such connection. As before, we model the maximum bid from other bidders (i.e. not the Agent) by introducing 
an i.i.d. family of nonnegative random variables $(B^{\nta}_k)_{k\in\N}$, where  $B^{\nta}_k$ represents the maximal bid of other bidders during the $k$-non-targeted advertising auction (in all the population).  The paying rule on this non-targeted auction is defined by a function $\bc^{\nta}:\R^2\rightarrow \R$, which is also assumed to be either of first-price or second-price auction rule. 
We stress that the auction rules used for the targeted advertising and the non-targeted advertising auctions do not necessarily have to be the same.
\end{itemize}

\paragraph{Advertising bidding map strategies.}

Given that there are now $M$ individuals, targeted advertising, and non-targeted advertising, a general bidding map control will take a more complex form with respect to  previous model. 
Informally, a bidding map control  is a random process, depending only upon past events (i.e. non anticipative), and valued in $\R^{M+1}$. The idea is that this vector process 
will store the $M$ bids that the Agent would like to make for each individual $m\in \llbracket 1,M\rrbracket$ if he were to connect to a website displaying targeted ads, and the remaining coordinate corresponds to the bid that the Agent would like to make if someone (anonymous) connects to a website using non-targeted advertising. Therefore,  $M+1$ potential bids are required at any time, hence the term {\it bidding map}.
    

To formalise this concept,  let us introduce the filtration $\F=(\Fc_t)_{t\in \R_+}$ generated by the processes 
\beqs 
((N^{m,\i}, N^{m,\ta}, N^{m,\nta}, N^{m,\danger}, B^{m,\ta}_{N^{m,\ta}}, N^{m,S})_{m\in \llbracket 1,M\rrbracket}, B^{\nta}_{N^{\nta}}, ((N^{m,i,\soc})_{m,i\in \llbracket 1,M\rrbracket})
\enqs 
where $N^{\nta}:=\sum_{m=1}^M N^{m,\nta}$ globally counts the connections to a website displaying non-targeted ads.
An open-loop bidding map control, denoted by $\beta$ $\in$ $\Pi_{OL}$,   is  a process $\beta=(\beta^m)_{m=0,...,M}$, valued in $\R^{M+1}_+$, predictable and progressively measurable w.r.t. the filtration $\F$. 
For $m=1,...,M$, $\beta^m_t$ is the bid that the Agent would make if a targeted advertising auction for individual $m$ happened at time $t$. The remaining coordinate, $\beta^0_t$ is the bid that the Agent would make if a non-targeted advertising auction occurs  at time $t$. In other words,  if an individual connects to a website displaying targeted ads (resp. non-targeted ads), the website will open the targeted advertising auction for this individual (resp. the non-targeted auction for this connection), look at the bidding map $\beta_t=(\beta^m_t)_{m=0,..., M}$, and automatically use the bid $\beta^m_t$ for the connected individual $m$ $\in$ 
$\llbracket 1,M\rrbracket$ (resp.  the bid $\beta^0_t$)  inscribed in this bidding map as the bid of the Agent for this auction. This allows the agent to specify a different bid for each individual, which encodes the idea of {\it targeted}-advertising, or a bid that do not depend upon who  is connecting, which encodes the idea of {\it non-targeted} advertising. 

\paragraph{The information dynamic.}

Given an open-loop bidding map control $\beta$, we define the information dynamic process  $X^{m,\beta}$ valued in $\{0,1\}$ of  any individual $m$ $\in \llbracket 1,M\rrbracket$ of the population as follows: 
$$
\begin{cases}
X^{m,\beta}_{0^-} &=\;  0,  \\
dX^{m,\beta}_t & = \; (1-X^{m,\beta}_{t-})(dN^{m,\i}_t+{\bf 1}_{\beta^m_t\geq B^{m,\ta}_{N_t^{m,\ta}}}dN^{m,\ta}_t \\
  & \quad \quad \quad \quad \quad + \; {\bf 1}_{\beta^0_t\geq B^{\nta}_{N_t^{\nta_t}}}dN^{m,\nta}_t+\sum_{i=1}^M X^{i,\beta}_{t-}dN^{m,i,\soc}_t), \quad t \geq 0. 
\end{cases}
$$
The interpretation of this dynamic is similar to previous sections for the first two terms (but they are now related to a given individual  $m\in \llbracket 1,M\rrbracket$), and  with additional terms which are 
essentially related to the new features of non-targeted advertising and social interactions.
Each individual $m$ starts uninformed ($X^{m,\beta}_{0^-}=0$). Once individual $m$ is informed ($X^{m,\beta}_t=1$), he stays informed ($(1-X^{m,\beta}_{t-})$ part). As long as he is not informed,  he can receive the information either
by connecting to  to a website intrinsically containing the Information ($dN^{m,\i}_t$ part), or by connecting to a website displaying targeted ads ($dN^{m,\ta}_t$ part), and when the Agent's ad is displayed to him, i.e. iff the Agent wins the targeted advertising auction (${\bf 1}_{\beta^m_t\geq B^{m,\ta}_{N_t^{m,\ta}}}$ part). 
Furthermore,  individual $m$  has the possibility  to 
\begin{itemize}
    \item  browse through websites displaying non-targeted ads ($dN^{m,\nta}_t$ part), in which case he will get informed if and only if the Agent's ad is displayed to him, i.e. iff the Agent wins the non-targeted advertising auction (${\bf 1}_{\beta^0_t\geq B^{\nta}_{N_t^{\nta_t}}}$ part), 
    \item  and socially interact with individual $i$ ($dN^{m,i,\soc}_t$ part). In this case, he will get informed whenever individual $i$ is informed ($X^{i,\beta}_t$ part). 
\end{itemize}

\paragraph{Proportion-based bidding policy.} A proportion-based bidding policy is defined by a pair of functions $\mathfrak{b}=(\mathfrak{b}^{\ta},\mathfrak{b}^{\nta})$ defined both from 
$\P_M$ $:=$ $\{ \frac{k}{M}: k=0,\ldots,M-1\}$ 
into $\R_+$. 
To any such policy we associate the open-loop bidding map control $\beta^{\mathfrak{b}}$ satisfying the feedback form  constraint 
\beqs 
\beta^{m,\mathfrak{b}}_t&=&\mathfrak{b}^{\ta}\Big(\frac{1}{M}\sum_{i=1}^M X^{i,\beta^{\mathfrak{b}}}_{t-}\Big)(1-X^{m,\beta^{\mathfrak{b}}}_{t-}), \quad m\in \llbracket 1,M\rrbracket,\\
\beta^{0,\mathfrak{b}}_t&=&\mathfrak{b}^{\nta}\Big(\frac{1}{M}\sum_{i=1}^M X^{i,\beta^{\mathfrak{b}}}_{t-}\Big){\bf 1}_{\frac{1}{M}\sum_{i=1}^M X^{i,\beta^{\mathfrak{b}}}_{t-}<1}, \quad \quad  t \geq 0.  
\enqs 
In other words, a bidding map control associated to a proportion-based bidding policy  formalises a strategy where in the targeted auction, the Agent  makes a bid 
for an individual $m$ that depends only on the proportion of informed people $\frac{1}{M}\sum_{i=1}^M X^{i,\beta^{\mathfrak{b}}}_{t-}$ at any time $t$,  and whether the individual $m$ is informed or not, and where in the non-targeted auction, she makes a bid depending only on  the proportion of informed people.

\paragraph{Cost function.} 
Given a bidding map control $\beta$, the expected cost incurred to the Agent is  defined by
\begin{align}  
V(\beta) &= \;  \E\Big[\sum_{m=1}^M\Big(\int_0^\infty K(1-X^{m,\beta}_{t-})dN^{m,\danger}_t+\int_0^\infty {\bf 1}_{\beta^m_t\geq B^{m,\ta}_{N_t^{m,\ta}}} \bc^{\ta}(b^{m,\ta}_t,B^{m,\ta}_{N_t^{m,\ta}}) d N_t^{m,\ta}  \\
& \hspace{3cm} + \; \int_0^\infty {\bf 1}_{\beta^0_t\geq B^{\nta}_{N^{\nta}_t}} \bc^{\nta} (\beta^0_t,B^{\nta}_{N^{\nta}_t}) dN^{m,\nta}_t\Big)\Big].  \label{defVsocialno} 
\end{align}  
This cost function is similar to previous model in Section \ref{secsocialdiscount}  except that there is a cost for each individual $m\in \llbracket 1,M\rrbracket$ in the population, ($\sum_{m=1}^M$ part), and that there is an additional term 
$\int_0^\infty {\bf 1}_{\beta^0_t\geq B^{\nta}_{N^{\nta}_t}} \bc^{\nta}(\beta^0_t,B^{\nta}_{N^{\nta}_t}) dN^{m,\nta}_t$ that measures the the non-targeted advertising cost of the strategy. 

\vspace{3mm}

We now state the main  result for this model.

\begin{Theorem}\label{theo-social-no-discount}
The minimal  cost is given by
 \beqs 
V^\star \; := \; \inf_{\beta\in \Pi_{OL}}V(\beta) &=& 
\sum_{p\in\P_M} 
v(p),
\enqs 
where 
$v(p)=\inf_{b^{\ta},b^{\nta}\in\R_+} v^{b^{\ta},b^{\nta}}(p)$, 
with 
\beqs 
v^{b^{\ta},b^{\nta}}(p)&=& \frac{K+  \eta^{\ta}\E\big[\bc^{\ta}(b^{\ta},B_1^{\ta}){\bf 1}_{b^{\ta}\geq B^{1,\ta}_1}\big]+\eta^{\nta}\E\Big[\frac{\bc^{\nta}(b^{\nta},B_1^{\nta})}{1-p}{\bf 1}_{b^{\nta}\geq\frac{B^{1,\nta}_1}{1-p}}\Big] }{\eta^{\i}+\eta^{\ta}\P\big[b^{\ta}\geq B^{1,\ta}_1\big] 
+\eta^{\nta}\P\big[b^{\nta}\geq B^{\nta}_1\big]+p\eta^{\soc}}. 
\enqs 
For all $p\in\P_M$, 
the set 
$\argmin_{b^{\ta},b^{\nta}\in\R_+} v^{b^{\ta},b^{\nta}}(p)$ 
is not empty, and any
proportion-based  bidding policy $\mfb^\star$ $=$ $(\mfb^{\star,\ta},\mfb^{\star,\nta})$ such that
\begin{align} \label{bstarinf} 
(\mfb^{\star,\ta}(p),\mfb^{\star,\nta}(p)) &\in \;  
\argmin_{b^{\ta},b^{\nta}\in\R_+} 
v^{b^{\ta},b^{\nta}}(p), \quad \forall 
p \in  \P_M. 
\end{align} 
yields an optimal bidding map control $\beta^{\mfb^\star}$. Moreover, in the case with fully second-price auctions, i.e. $\bc^{\ta}(b,B)=\bc^{\nta}(b,B)=B$, we have 
$v(p)=\inf_{b\in\R_+} v^{b,(1-p)b}(p)$. 
The set 
$\argmin_{b\in\R_+} v^{b,(1-p)b}(p)$
 is not empty, and any 
proportion-based  bidding policy defined by 
%
$\mfb^{\star,\ta}(p)$ $=$ $\mrb^\star(p)$, $\mfb^{\star,\nta}(p)$ $=$ $(1-p)\mrb^{\star}(p)$ with 
\beqs 
\mrb^\star(p) &=& 
\argmin_{b\in\R_+} 
v^{b,(1-p)b}(p), \quad \forall p\in \P_M, 
\enqs 
yields an optimal bidding map control $\beta^{\mfb^\star}$.
\end{Theorem}

\noindent{\bf Interpretation.} Let us provide some  interpretations of the formulas in Theorem \ref{theo-social-no-discount}. 
\begin{itemize}
    \item {\it The sum part ``$\sum_{p\in\P_M}$". 
    } 
    We can split the problem in several successive problems each consisting in optimally going from a proportion $\frac{k}{M}$ of informed people to a proportion  $\frac{k+1}{M}$, for $k\in\{0,...,M-2\}$. The fact that there is no discount rate implies that the time when each problem starts does not matter, and therefore,  these successive problems can be optimised independently, i.e. one by one. 
    \item {\it The term in the sum.} The justification of the form of the terms in the sum is similar to the justification given for the previous model: the fraction can be split into two fractions, one corresponding to the expected cost perceived during this period, and the other one corresponding to the expected cost perceived at the termination time of this period.
    \item {\it The term $\frac{\bc^{\nta}(b,B_1^{\nta})}{1-p}$.} Notice that in the formula, $B^{1,\ta}_1$ and $\frac{B^{\nta}_1}{1-p}$ play symmetric roles. It is {\it as if} the non-targeted advertising mechanism with price $B^{\nta}_1$ is  equivalent to a targeted advertising mechanism with price $\frac{B^{\ta}_1}{1-p}$. In other words, making the advertising mechanism not targeted essentially is equivalent to multiply the ad cost by $\frac{1}{1-p}$. This is natural since when the ad mechanism is not targeted, there is a probability $p$ that it displays the ad to an already informed individual. This means that only a proportion $1-p$ of the paid ads will be useful, and thus, for each useful ad, an average number of $\frac{1}{1-p}$ ads has to be paid (including the useful one).  In other words, we have to pay the price of $\frac{1}{1-p}$ ads to display an ad to an uninformed individual.
    \item {\it The term $p\eta^{\soc}$.} Notice that in the formula, $p\eta^{\soc}$ plays the same role as $\eta^{\i}$. This is consistent with the intuition that socially interacting with an informed individual has the same effect as visiting a website containing the information: it will inform the individual and not cost anything to the Agent. The more individuals are informed, the more likely such interaction is to occur. More precisely, each informed individual ``plays the role'' of a website containing the information, such that an individual has intensity $\frac{1}{M}\eta^{\soc}$ to ``visit'' it, and thus, with a $k$ informed individuals, it yields an intensity $\frac{k}{M}\eta^{\soc}=p\eta^{\soc}$.
\end{itemize}

\vspace{1mm}

We introduce the following special optimal proportion-based  bidding policy.

\begin{Definition}[Smallest optimal proportion-based  bidding policy]\label{def-smallest-proportion-based-bid}
There exists a unique  proportion-based bidding policy $\mfb_{min}^\star$ $=$ $(\mfb_{min}^{\star,\ta},\mfb_{min}^{\star,\nta})$   such that any proportion-based bidding policy $\mfb^{\star}$ $=$ $(\mfb^{\star,\ta},\mfb^{\star,\nta})$  
as in Theorem \ref{theo-social-no-discount} satisfies $\mfb^{\star,\ta}_{min}(p)$ $\leq$ $\mfb^{\star,\ta}(p)$ and $\mfb^{\star,\nta}_{min}(p)$ $\leq$ $\mfb^{\star,\nta}(p)$ for all $p\in\P_M$. 
$\mfb_{min}^{\star}$ is called the smallest optimal proportion-based  bidding policy.
\end{Definition}

\begin{Remark}
The above result comes from the identity
\beqs
& & \argmin_{b^{\ta},b^{\nta}\in\R_+} v^{b^{\ta},b^{\nta}}(p) \\
&=& \argmax_{b^{\ta}\in \R_+}\E\Big[(v(p)-B^{1,\ta}_{1}){\bf 1}_{b^{\ta}\geq B^{1,\ta}_{1}}\Big] \times \argmax_{b^{\nta}\in \R_+}\E\Big[\big((1-p)v(p)-B^{\nta}_{1}\big){\bf 1}_{b^{\nta}\geq B^{\nta}_{1}}\Big]
\enqs
which follows  from a Bellman and verification result type property, essentially allowing to see, in this dynamic problem, the long term optimal bid for a targeted advertising auction (resp. for a non-targeted advertising auction) as a greedy optimal bid for a static auction with immediate reward $v(p)$ (resp. $(1-p)v(p)$) when the auction is won. Consequently, we have, for all $p\in\P_M$, 
\beqs 
\mfb^{\star,\ta}_{min}(p)&=&\min \argmax_{b^{\ta}\in \R_+}\E\Big[(v(p)-B^{1,\ta}_{1}){\bf 1}_{b^{\ta}\geq B^{1,\ta}_{1}}\Big],\\
\mfb^{\star,\nta}_{min}(p)&=&\min \argmax_{b^{\nta}\in \R_+}\E\Big[\big((1-p)v(p)-B^{\nta}_{1}\big){\bf 1}_{b^{\nta}\geq B^{\nta}_{1}}\Big].
\enqs 
It is also possible to see, from the proof of Theorem \ref{theo-social-no-discount}, that the open-loop bidding map control $\beta^{\mfb_{min}^\star}$ is the smallest optimal open-loop bidding map control, i.e. for all optimal open-loop bidding map control 
$\beta$ $=$ $(\beta^m)_{m\in\llbracket 0,M\rrbracket}$, we have $\beta^{m,\mfb_{min}^\star}_t$ $\leq$ $\beta^m_t$ for all $m\in \llbracket 0, M\rrbracket$, $(\omega,t)$-a.e., i.e. for almost every $(\omega,t)\in \Omega\times \R_+$ w.r.t. the measure $\P\otimes \Bc(\R_+)$.
\end{Remark}

We have the following properties  about the sensitivity to parameters and upper bounds of the optimal value and smallest minimal optimal bidding policy.
\begin{Proposition}\label{propsocialno-sensitivity}
The optimal value $V^\star$  and the smallest optimal proportion-based  bidding policies $b^{\star,\ta}_{min}(p)$ and $b^{\star,\nta}_{min}(p)$ for targeted and for non-targeted advertising are decreasing in $\eta^{\i},\eta^{\ta}, \eta^{\nta}$, and $\eta^{\soc}$. Furthermore,  $b^{\star,\nta}_{min}(p)$  is decreasing in $p$, while 
$b^{\star,\ta}_{min}(p)$ is 
    \begin{itemize}
        \item decreasing in $p$ when  there is no non-targeted advertising ($\eta^{\nta}=0$),
        \item increasing in $p$ when  there is no social interactions ($\eta^{\soc}=0$).
    \end{itemize}
Finally, we have for all $p$ $\in$ $[0,1]$, 
\beqs 
b^{\star,\ta}_{min}(p), \; b^{\star,\nta}_{min}(p)   \;  \leq \;  v(p) \; \leq \; \frac{K}{\eta^{\i}+p\eta^{\soc}}. 
\enqs 
\end{Proposition}

\vspace{1mm}

The interpretation of the above proposition related to the monotonicity of the optimal value and smallest proportion-based bidding policies with respect to the intensity parameters is similar as in the previous section. Let us discuss the monotonicity properties with respect to the proportion $p$ of informed individuals.  Recall that the agent has at disposal four channels of information: (1) her website which informs without cost the individual, (2) the social interaction which also informs with probability $p$ the individual and without cost, (3) the targeted ad which informs surely the individual with cost $\bc^{\ta}$, and (4) the non-targeted ad which informs with probability $1-p$ the individual with cost $\bc^{\nta}$.  It is then a trade-off to choose the most efficient channel for diffusing the information. So, in the case where the fourth channel is not accessible ($\eta^{\nta}$ $=$ $0$), and when $p$ is increasing, this will only affect positively the efficiency of the social interactions (second channel), and therefore the agent will bid less for the targeted ad. On the other hand, in the case where the second channel is absent ($\eta^{\soc}$ $=$ $0$), and when $p$ is increasing, this will only affect the fourth channel, which loses in efficiency since many individuals are already informed. Consequently, it becomes more interesting to bid on the targeted ad. In the general case, when $p$ is increasing, the second channel gains in efficiency but the fourth channel  become more costly.

\begin{Remark} The monotonicity properties of the smallest optimal proportion-based bidding policies with respect to $p$ has useful implications regarding their  computational cost. Indeed, the practical implementation of these optimal bids require to compute (via the search of the infimum in \eqref{bstarinf}) 
$(\mfb_{min}^{\star,\ta}(p),\mfb_{min}^{\star,\nta}(p))$ for any $p$ $\in$ $\P_M$. This is a priori very  expensive for large $M$. However, by taking advantage of the monotonicity in $p$ of $(\mfb_{min}^{\star,\ta}(p),\mfb_{min}^{\star,\nta}(p))$, one can  considerably reduce the computational complexity.  For  instance, if one starts by computing optimal bids for $p= \frac{1}{2}$, the computation of $\mfb^{\star,\nta}_{min}(\frac{1}{2})$ allows to limit the search for $\mfb^{\star,\nta}_{min}(p)$ to $[0,\mfb^{\star,\nta}_{min}(\frac{1}{2})]$ for $p>\frac{1}{2}$, and to 
$[\mfb^{\star,\nta}_{min}(\frac{1}{2}), \frac{K}{\eta^{\i}+p\eta^{\soc}}]$ for $p<\frac{1}{2}$. In particular, one can search $\mfb^{\star,\nta}_{min}(\frac{3}{4})$ in $[0,\mfb^{\star,\nta}_{min}(\frac{1}{2})]$ and $\mfb^{\star,\nta}_{min}(\frac{1}{4})$ in $[\mfb^{\star,\nta}_{min}(\frac{1}{2}), \frac{K}{\eta^{\i}+\frac{1}{4}\eta^{\soc}}]$. 
The computation of $(\mfb^{\star,\nta}_{min}(p))_{p\in\frac{\llbracket 0, M\llbracket}{M}}$ can thus clearly be made by dichotomy. For example, let us assume that there is only non-targeted advertising, and to simplify, let us assume that $M=2^L$ for some $L\in\N$. Then, only $L=\log_2(M)$ dichotomies have to be made, and at the $\ell$-th dichotomy, there are $2^\ell$ minimizers to find in $2^\ell$ consecutive intervals with total (i.e. cumulative) length upper bounded by $\frac{K}{\eta^{\i}}$. Assuming that the computational time of the search for a minimizer is proportional to the length of the interval on which it is searched, the computational complexity of each dichotomy iteration is thus $\Oc(\frac{K}{\eta^{\i}})$, and therefore, the total computational complexity of the whole minimal bidding policy is $\Oc(\frac{K}{\eta^{\i}}\log_2(M))$ and is thus only logarithmic in the size $M$ of the population, which suggests that this algorithm should be tractable even for large population.
\end{Remark}

\begin{Remark}[Mean-field approximation] \label{rem-meanfield} 
As in any population models with enough symmetry, it is expected that when $M$ gets large, we obtain a {\it mean-field limit}. Let us check  formally how  it is derived.  
Notice that the mean over the population of the Agent's optimal value is equal to  
 \beqs 
\frac{1}{M} V^\star  &=&  \frac{1}{M}\sum_{p\in \P_M} v(p),
\enqs 
which thus takes the form of a Riemann sum, hence when $M$ $\rightarrow$ $\infty$, 
\beqs 
\frac{1}{M} V^\star &\simeq &  \int_0^1 v(p)dp. 
\enqs 
where $v$ is extended on $[0,1)$ with the same expression as in Theorem \ref{theo-social-no-discount}. Such result can be useful for two reasons:
\begin{enumerate}
    \item To obtain an analytical approximation of the optimal value in some cases where the integral can be explicitly computed,
    \item and to provide a way to numerically approximate the optimal value, by discretising the integral with a suitable discretisation step. This can be useful with very large population, where one might want to speed up the computation.
\end{enumerate}
It is also possible to formally derive a differential optimal control problem on the proportion of informed users $(p_t)_{t\in\R_+}$ such that the optimal value and optimal control from this model are the limit of the corresponding objects in our model when $M\rightarrow \infty$. The controlled dynamic is defined by
\beqs 
\frac{dp^\beta_t}{dt}= (1-p^\beta_t)\big(\eta^{\i}+ \eta^{\ta}\P[\beta^{\ta}_t\geq B^{\ta}_1] + \eta^{\nta}\P[\beta^{\nta}_t\geq B^{\nta}_1] +\eta^{\soc} p^\beta_t\big),
\enqs 
with deterministic control $\beta=(\beta^{\ta},\beta^{\nta})$,  and  with  cost functional
\beqs 
V(\beta) &=& \int_0^{\infty}\Big\{ (1-p^\beta_t)\big(K+\eta^{\ta}\E[\bc^{\ta}(\beta^{\ta}_t,B^{\ta}_1){\bf 1}_{\beta^{\ta}_t\geq B^{\ta}_1}]\big) \\
& & \hspace{2cm} + \; \eta^{\nta}\E\big[\bc^{\nta}(\beta^{\nta}_t,B^{\nta}_1){\bf 1}_{\beta^{\nta}_t\geq B^{\nta}_1}\big] \Big\}dt.
\enqs 
\end{Remark}





\section{Examples with explicit computations} \label{sec:example}

Notice that the results of all the models presented in this work have a solution (optimal value and policy) that can be expressed in the form
\beqs 
\inf_{b\in\R_+}/\sup_{b\in\R_+}/\argmin_{b\in\R_+}/\argmax_{b\in\R_+}\frac{a_1+a_2\E[\bc(b,B^{\ta}_1){\bf 1}_{b\geq B^{\ta}_1}]+a_3\E[c(\bb,B^{\nta}_1){\bf 1}_{\bb\geq B^{\nta}_1}]}{a_4+a_5\P[b\geq B^{\ta}_1] + a_6\P[\bb\geq B^{\nta}_1]}
\enqs 
with well chosen parameters $(a_i)_{i\leq 6}$. In this section, we  discuss two types of distributions for $B^{\ta}_1$ and $B^{\nta}_1$ that will lead to fully explicit formulas for the optimal bidding policy.
\subsection{Constant maximal bid from other bidders}
We consider the case where the maximal bids from other bidders, i.e. $(B^{\ta}_k)_{k\in\N}$ for the targeted advertising auctions, and $(B^{\nta}_k)_{k\in\N}$ for the non-targeted advertising auctions, are constant, i.e. $B^{\ta}_k=B^{\ta}\in\R_+$ and $B^{\nta}_k=B^{\nta}\in\R_+$. Under this assumption, the first-rice auction or second-price auction cases essentially become equivalent, and we focus on the second price type of auction, i.e. the auction payment rule $\bc(b,B)=B$. Let us study two cases:
\begin{enumerate}
    \item The commercial advertising problem with purchase-based gain function, and 
    \item The social marketing problem with no discount factor and with social interactions and non-targeted advertising.
\end{enumerate}

\subsubsection{Commercial advertising  with purchase-based gain function}
In this case, we have
\beqs 
V(\beta^b)&=& \frac{\eta^{\i}K+\eta^{\ta}(K-B^{\ta}) {\bf 1}_{b\geq B^{\ta}}}{\eta^{\i}+ \rho + \eta^{\ta}{\bf 1}_{b\geq B^{\ta}}},
\enqs 
and any $b^\star$ $\in$ $\argmax_{b\in\R_+}$ $V(\beta^b)$ yields an optimal constant bid. Notice that $V(\beta^b)$ only takes two possible values, one for $b<B^{\ta}$ and one for $b\geq B^{\ta}$. The optimisation thus reduces to choose either $b<B^{\ta}$ (for instance $b=0$), either $b\geq B^{\ta}$ (for instance $b=B^{\ta}$). 
Hence, an optimal bid is $b^\star\geq B^{\ta}$ iff
\beqs 
\frac{\eta^{\i}K+\eta^{\ta}(K-B^{\ta}) }{\eta^{\i}+ \rho + \eta^{\ta}}> \frac{\eta^{\i}K}{\eta^{\i}+ \rho },
\enqs 
which can be rewritten equivalently after some straightforward calculation as
\beqs 
B^{\ta} \; < \; \frac{\rho}{\eta^{\i}+ \rho}K.
\enqs 
We clearly see what we had already established in the general case: the optimal bids are ``decreasing'' in $\eta^{\i}$ and ``increasing'' in $\rho$. Namely,  the smallest optimal bid is $B^{\ta}{\bf 1}_{B^{\ta}\leq \frac{\rho}{\eta^{\i}+ \rho}}$, which is  clearly a decreasing (resp. increasing) function of $\eta^{\i}$ (resp. of  $\rho$). 

There is another interesting optimal bid, that is, the bid  $\frac{\rho}{\eta^{\i}+ \rho}K$. Indeed, this bid is the only one to be optimal {\it regardless} $B^{\ta}$. In other words, by  assuming (or knowing) that other bidders' maximal bid is constant, 
we obtain  a dominant bidding strategy $\frac{\rho}{\eta^{\i}+ \rho}K$, which is optimal whatever (hence robust to) the constant value of $B^{\ta}$.

\subsubsection{Social marketing problem with no discount factor and with social interactions and non-targeted advertising}

In this case, the optimal bidding map control is  obtained from a proportion-based bidding policy $\mfb^\star$ $=$ $(\mfb^{\star,\ta}, \mfb^{\star,\nta})$ with $\mfb^{\star,\ta}(p)$ $=$  $\mrb^\star(p)$ and $\mfb^{\star,\nta}(p)$ $=$ $(1-p)\mrb^\star(p)$, where
\beqs 
\mrb^\star(p) & \in & \argmin_{b\in \R_+}\frac{K+\eta^{\ta}B^{\ta}{\bf 1}_{b\geq B^{\ta}}  +\eta^{\nta}\frac{B^{\nta}}{1-p}{\bf 1}_{b\geq \frac{B^{\nta}}{1-p}} }{\eta^{\i}+\eta^{\ta}{\bf 1}_{b\geq B^{\ta}}+\eta^{\nta}{\bf 1}_{b\geq \frac{B^{\nta}}{1-p}} +p\eta^{\soc}}. 
\enqs 
In order to obtain simple and interpretable formula,   let us assume that there is only one type of advertising.

\paragraph{Only targeted advertising.}  If there is only targeted advertising, i.e. if $\eta^{\nta}=0$, we have
\beqs 
\mfb^{\star,\ta}(p) & \in & \argmin_{b\in \R_+} \frac{K+\eta^{\ta}B^{\ta}{\bf 1}_{b\geq B^{\ta}}}{\eta^{\i}+\eta^{\ta}{\bf 1}_{b\geq B^{\ta}} +p\eta^{\soc}}. 
\enqs 
Here again, we are reduced to compare two costs: 
\beqs 
\frac{K}{\eta^{\i}+p\eta^{\soc}} & \text{ and } & \frac{K+\eta^{\ta}B^{\ta}}{\eta^{\i}+\eta^{\ta} +p\eta^{\soc}},
\enqs 
the first one being obtained for $b<B^{\ta}$, and the second one for $b\geq B^{\ta}$. The best option will be $\mfb^{\star,\ta}(p)$ $\geq$ $B^{\ta}$ if and only if
\beqs 
\frac{K}{\eta^{\i}+p\eta^{\soc}} &>&  \frac{K+\eta^{\ta}B^{\ta}}{\eta^{\i}+\eta^{\ta} +p\eta^{\soc}},
\enqs 
which is equivalent to
\beqs 
p \; < \; \frac{\frac{K}{B^{\ta}}-\eta^{\i}}{\eta^{\soc}} \; =: \; p_\star. 
\enqs 
This means that before the proportion of informed individuals attains the threshold  informed proportion $p_\star$, one should bid higher than $B^{\ta}$ (and thus display ads), and when the informed proportion  exceeds $p_\star$, one should bid lower than $B^{\ta}$ (and thus stop displaying ads). 
Assuming that $K$ $\geq$ $\eta^{\i}B^{\ta}$, we also notice  that this threshold $p_\star$ is decreasing in $\eta^{\i}$ and in $\eta^{\soc}$. This is interpreted as follows: 
\begin{itemize}
    \item First of all, the fact that there is an informed proportion threshold under which the Agent should display ads and above which she should not display ads necessarily comes from the social interactions. 
    Indeed, in  no-social interaction case ($\eta^{\soc}=0$), this threshold $p_\star$ $=$ $\infty$. 
    The fact that the presence of social interactions is susceptible to introduce a finite threshold proportion above  which the Agent should stop displaying ads, can be also understood as follows: suppose that an individual connects to a website displaying targeted ads. With social interactions, the more people are informed, the sooner this individual will learn the Information anyway, by interacting with an informed individual. Therefore, the incentive of the Agent to display the ad to him is weaker as the proportion of informed individuals increases, which justifies that the bid she is willing to make is smaller, and once it is small enough to fall below $B^{\ta}$, the Agent will stop displaying ads.
    \item The interpretation of the decreasing nature of the threshold proportion in $\eta^{\i}$ and $\eta^{\soc}$ is the following. For a fixed proportion of informed individuals $p$, increasing the intensity of social interactions $\eta^{\soc}$ will also make more likely a soon interaction with an informed people, thus weakening the Agent's incentive to display an ad, such that this incentive will be fully compensated after a smaller informed proportion. Likewise, increasing the intensity $\eta^{\i}$ of connections to a website containing the information will make people inform themselves faster, thus catalyzing the increase of the informed proportion, in turn decreasing the Agent's incentive to display an ad.
\end{itemize}

\vspace{2mm}

Finally, we can estimate  the optimal value when $M$ $\rightarrow$ $\infty$,  by using the mean-field approximation as observed in Remark \ref{rem-meanfield}. First, consider the case when  $p_\star$ $\geq$ $1$. In this case, we have  $\mfb^{\star,\ta}(p)$ $\geq$ $B^{\ta}$, for all $p$ $\in$ $[0,1]$, and so 
\beqs 
\frac{1}{M}V^\star & \simeq & \int_0^1 v(p) dp \; = \;  \int_0^1    \frac{K+\eta^{\ta}B^{\ta}{\bf 1}_{\mfb^{\star,\ta}(p)\geq B^{\ta}}}{\eta^{\i}+\eta^{\ta}{\bf 1}_{\mfb^{\star,\ta}(p)\geq B^{\ta}} +p\eta^{\soc}} dp  
\; = \; \int_0^{1}\frac{K+\eta^{\ta}B^{\ta}}{\eta^{\i}+\eta^{\ta} +p\eta^{\soc}} dp \\
&=&\frac{K+\eta^{\ta}B^{\ta}}{\eta^{\soc}}\ln\Big[\frac{\eta^{\i}+\eta^{\ta} +\eta^{\soc}}{\eta^{\i}+\eta^{\ta}}\Big]. 
\enqs 
When $0 \leq p_\star < 1$, we have 
\beqs 
\frac{1}{M} V^\star  &\simeq &  \int_0^1    \frac{K+\eta^{\ta}B^{\ta}{\bf 1}_{\mfb^{\star,\ta}(p)\geq B^{\ta}}}{\eta^{\i}+\eta^{\ta}{\bf 1}_{\mfb^{\star,\ta}(p)\geq B^{\ta}} +p\eta^{\soc}} dp   
\; = \;  \int_0^{p_\star}\frac{K+\eta^{\ta}B^{\ta}}{\eta^{\i}+\eta^{\ta} +p\eta^{\soc}} dp \; + \; \int_{p_\star}^1\frac{K}{\eta^{\i}+p\eta^{\soc}}dp\\
&=& \frac{K+\eta^{\ta}B^{\ta}}{\eta^{\soc}}\ln\Big[\frac{\eta^{\i}+\eta^{\ta} +p_\star\eta^{\soc}}{\eta^{\i}+\eta^{\ta}}\Big] \; - \; \frac{K}{\eta^{\soc}}\ln\Big[ \frac{\eta^{\i} +p_\star\eta^{\soc}}{\eta^{\i}+\eta^{\soc}}\Big] \\
&=& \frac{K+\eta^{\ta}B^{\ta}}{\eta^{\soc}}\ln\Big[ \frac{\eta^{\ta} +\frac{K}{B^{\ta}}}{\eta^{\i}+\eta^{\ta}}\Big] \; - \; \frac{K}{\eta^{\soc}}\ln\Big[ \frac{\frac{K}{B^{\ta}}}{\eta^{\i}+\eta^{\soc}}\Big].  \\
\enqs

\paragraph{Only non-targeted advertising.} If there is only non-targeted advertising, i.e. if $\eta^{\ta}=0$, we have
\beqs 
\mfb^{\star,\nta}(p) &\in & \argmin_{b\in \R_+}\frac{K  +\eta^{\nta}\frac{B^{\nta}}{1-p}{\bf 1}_{b\geq B^{\nta}} }{\eta^{\i}+\eta^{\nta}{\bf 1}_{b\geq B^{\nta}
} +p\eta^{\soc}}. 
\enqs 
Here again, we are reduced to compare two costs: 
\beqs 
\frac{K}{\eta^{\i}+p\eta^{\soc}} \;\;\; \text{ and } \;\;\;  \frac{K  +\eta^{\nta}\frac{B^{\nta}}{1-p}}{\eta^{\i}+\eta^{\nta} +p\eta^{\soc}},
\enqs 
the first one being obtained for $b<B^{\nta}$, and the second one for $b\geq B^{\nta}$. The best option will be $\mfb^{\star,\nta}(p)$ $\geq$ $B^{\nta}$ if and only if
\beqs 
\frac{K}{\eta^{\i}+p\eta^{\soc}} &>& \frac{K  +\eta^{\nta}\frac{B^{\nta}}{1-p}}{\eta^{\i}+\eta^{\nta} +p\eta^{\soc}},
\enqs 
which is equivalent to
\beqs 
p &<& \frac{K-\eta^{\i} B^{\nta}}{K+\eta^{\soc} B^{\nta}} \; =: \; \bar p_\star. 
\enqs 
This means that below  the informed proportion threshold $\bar p_\star$,  one should bid higher than $B^{\nta}$ (and thus display ads), and above  the informed proportion $\bar p_\star$, one should bid lower than $B^{\nta}$ (and thus stop displaying ads). 
When $K$ $\geq$ $\eta^{\i} B^{\nta}$, we notice, as in the ``only targeted advertising'' case, that  the informed proportion $\bar p_\star$  is decreasing in $\eta^{\i}$ and in $\eta^{\soc}$. 
The same interpretations  as in the ``only targeted advertising'' case  still hold, but there is an additional argument. Indeed, recall that in the ``only targeted advertising'' case, it is argued  that the presence of such threshold comes from the presence of social interactions, and that when they are absent ($\eta^{\soc}=0$), or more generally when $\eta^{\soc}$ is small enough, there is no threshold (the optimal bidding strategy is a constant bid). Here, notice that $\bar p_\star$ $<$ $1$, even  if $\eta^{\soc}=0$ (recall that we assumed that $\eta^{\i}>0$). Thus, as opposed to the previous ``only targeted advertising'' case, the existence of such threshold does not only come from social interactions. 
Displaying non-targeted ads always comes with the risk to display ads to already informed people, and thus paying for a useless ad. The more people are informed, the higher the risk. This explains why  after some proportion threshold, it is not worth to pay for displaying an ad, and thus the Agent has to stop doing so.

\subsection{Uniform maximal bid from other bidders}

We now consider the case where the other bidders' maximal bid is uniformly distributed.
Regarding the auction payment rule, we shall focus on the {\it first-price auction rule}, i.e., $\bc(b,B)=b$ (the same  argument applies to the {\it second-price auction rule}). 

We focus on the example of the purchase-based commercial advertising model, but the same argument can be adapted to the other models. From Theorem \ref{theo-purchase},  the gain value associated to a constant bid strategy is equal to 
\beqs 
V(\beta^b) &=&  
\frac{\eta^{\i}K+\eta^{\ta}(K-b)\P[b\geq B^{\ta}_1]}{\eta^{\i}+ \rho + \eta^{\ta}\P[b\geq B^{\ta}_1]}, \quad b \in \R_+.  
\enqs 
Denoting by  $[b^-,b^+]$ with  $b^-<b^+$, the support of the uniform distribution for  $B^{\ta}_1$, we can restrict the search for the argmax of $b$ $\mapsto$ $V(\beta^b)$ to  the interval $[b^-,b^+]$, and so 
\beqs 
b^\star & \in &  \argmax_{b\in [b^-,b^+]} \frac{\eta^{\i}K+\eta^{\ta}(K-b)\frac{b-b^-}{b^+-b^-}}{\eta^{\i}+ \rho + \eta^{\ta}\frac{b-b^-}{b^+-b^-}}. 
\enqs 
By making the change of variable
\beqs 
b' \; = \; \eta^{\i}+ \rho + \eta^{\ta}\frac{b-b^-}{b^+-b^-}, & \mbox{ i.e. } &  b \; = \;  \lambda_1 + \lambda_2 b'
\enqs 
with
\beqs 
\lambda_1 \; = \; b^--(b^+-b^-)\frac{\eta^{\i}+ \rho }{\eta^{\ta}},\quad \lambda_2\; = \; \frac{b^+-b^-}{\eta^{\ta}},
\enqs 
we see that $b^\star$ $=$ $\lambda_1 + \lambda_2 b^{',\star}$, with 
\begin{align} \label{b'} 
b^{',\star}
& \in  \; \argmax_{b'\in [b'^-,b'^+]} \frac{a_0+a_1b'+a_2b'^2}{b'} \; = : \;  \argmax_{b'\in [b'^-,b'^+]}  g(b'), 
\end{align} 
with
\beqs 
a_0 \; = \; \lambda_1(\eta^{\i}+\rho)-K\rho, \quad a_1 \; = \; K-\lambda_1+\lambda_2(\eta^{\i}+\rho), \quad a_2= -\lambda_2<0,
\enqs 
and 
\beqs 
b'^-=\eta^{\i}+ \rho, \quad b'^+=\eta^{\i}+ \rho+\eta^{\ta}. 
\enqs 
By writing the first-order condition for $g(b')$ in \eqref{b'}, we see that its derivative is equal to $a_2-\frac{a_0}{b'^2}$ which is negative, for $b'\in [b'^-,b'^+]\subset \R_+$, if and only if $b'^2\geq \frac{a_0}{a_2}$, and thus if and only if  $b'\geq \sqrt{\left(\frac{a_0}{a_2}\right)_+}$. 
The argmax for $g$ in $[b'^-,b'^+]$  is thus given by
\beqs 
b^{',\star} &=&  \max\Big[ b'^-, \min\Big(b'^+, \sqrt{\big(\frac{a_0}{a_2}\big)_+}\Big)\Big]
\enqs 
and thus
\beqs 
b^\star 
&=&  \max\big[b^-, \min(b^+, \bar b)\big]
\enqs 
where (after some straightforward calculation)
\beqs 
\bar b 
&=&  b^--(b^+-b^-)\frac{\eta^{\i}+ \rho }{\eta^{\ta}} +  \sqrt{\frac{b^+-b^-}{\eta^{\ta}}\Big(K\rho-b^-(\eta^{\i}+\rho)+\frac{b^+-b^-}{\eta^{\ta}}(\eta^{\i}+ \rho)^2\Big)_+}. 
\enqs

\section{Proof of main results} \label{sec:proofs}


We first prove the results  in the social marketing model with no discount rate,  and then show how the other results  can be reduced as particular cases of this model. 

\subsection{Proof of results in Section \ref{sec:nodiscount}}

\subsubsection{Proof of Theorem \ref{theo-social-no-discount}}  \label{sec:socialno}

Fix an arbitrary open-loop bidding map control $\beta$, and denote by 
\beqs 
p^{\beta}_t &=&  \frac{1}{M}\sum_{m=1}^M X^{m,\beta}_t,\quad t\in\R_+, 
\enqs 
the proportion of informed individuals at each time $t\in\R_+$.  The underlying idea of this proof is a change of variable from the numerous Poisson processes of the problem to the proportion $p^{\beta}_t$ in the cost function. 
In other words,  the suitable approach  to look at the problem is to express the optimisation of the cost $V(\beta)$  over the proportion $p^\beta$ running into $\{0,1/M,\ldots,1\}$ rather than over the times of jumps of the numerous Poisson processes $N^{\i}$, $N^{\danger}$, $N^{\ta}$, $N^{\nta}$  defined in our model.  
Notice that Poisson and proportion process are piece-wise constant processes, and thus the change of variable has to be done carefully.  To deal with this technical issue, we shall rely on  
the compensated processes of the Poisson processes, and use  martingale arguments in order to express $V(\beta)$ first with $dt$ thanks to the intensity processes, 
then make a continuous-time change of variable to obtain another intensity process, and then move back to jump processes with  $dp^\beta$. 



\vspace{1mm}

\noindent {\it Step 1: Intensity process of $p^\beta$.} From the dynamics of $X^\beta$, we have
\beqs 
dp^{\beta}_t 
&=& \frac{1}{M}\sum_{m=1}^M (1-X^{m,\beta}_{t-})\Big(dN^{m,\i}_t+{\bf 1}_{\beta^m_t\geq B^{m,\ta}_{N^{m,\ta}_t}}dN^{m,\ta}_t. \\
& & \hspace{3cm} + \; {\bf 1}_{\beta^0_t\geq B^{\nta}_{N^{\nta}_t}} dN^{m,\nta}_t+\sum_{i=1}^M  X^{i,\beta}_{t-}dN^{i,m,\soc}_t\Big) \\
&=&  \frac{1}{M}\sum_{m=1}^M (1-X^{m,\beta}_{t-})\Big(dN^{m,\i}_t+{\bf 1}_{\beta^m_t\geq B^{m,\ta}_{N^{m,\ta}_{t-}+1}}dN^{m,\ta}_t \\
& & \hspace{3cm} + \; {\bf 1}_{\beta^0_t\geq B^{\nta}_{N^{\nta}_{t-}+1}} dN^{m,\nta}_t+\sum_{i =1}^M  X^{i,\beta}_{t-}dN^{m,i,\soc}_t\Big). 
\enqs 
It follows that the process
\begin{align}\label{eq-martingale}
t &\longmapsto\;   p^{\beta}_t-\int_0^t\frac{1}{M}\sum_{m \in \llbracket 1,M \rrbracket}  (1-X^{m, \beta}_{s-})\Big(\;\eta^{\i}+{\bf 1}_{\beta^m_s\geq B^{m,\ta}_{N^{m,\ta}_{s-}+1}}\eta^{\ta}\\
&  \hspace{3cm}  + \; {\bf 1}_{\beta^0_s\geq B^{\nta}_{N^{\nta}_{s-}+1}} \eta^{\nta}+\sum_{i  \in \llbracket 1,M \rrbracket}  X^{i,\beta}_{s-}\eta^{\soc}\Big)ds, 
\end{align}
is a martingale. Indeed, a classical result of martingale theory for Point process,  is that for any Poisson process $N$ with intensity $\eta$, its  {\it compensated process}, defined by $(N_t-\eta t)_{t\in\R_+}$ is a martingale w.r.t. 
its natural  filtration, but also, clearly, w.r.t. any filtration generated by $N$ and any process $Y$ independent of $N$. This implies that all the Poisson processes considered in this model are martingales w.r.t. the filtration $\tilde{\F}=(\tilde{\Fc}_t)_{t\in\R_+}$ defined by
\beqs 
\tilde{\Fc}_t=\sigma((B^{m,\ta}_k)_{m\in\llbracket 1, M\rrbracket, k\in\N_\star},(B^{\nta}_k)_{k\in\N_\star}, (N^{m,\i}_s,N^{m,\ta}_s,N^{m,\nta}_s,N^{m,i,\soc}_s,N^{m,\danger}_s)_{m,i \in \llbracket 1,M \rrbracket, s\leq t}). 
\enqs 
Notice that  the processes in the integrand in \eqref{eq-martingale} is $\tilde{\F}$-predictable, which thus implies that the process in \eqref{eq-martingale} is a $\tilde{\F}$-martingale. Since $\Fc_t\subset \tilde{\Fc}_t$ for all $t\in\R_+$,  this implies that for any bounded positive $\F$-predictable process $H$, 
\beqs 
&& \E\Big[\int_0^\infty H_tdp^\beta_t\Big]\\
&=&\E\Big[\int_0^\infty H_t \frac{1}{M}\sum_{m \in \llbracket 1,M \rrbracket}   (1-X^{m, \beta}_t)\Big(\eta^{\i}+{\bf 1}_{\beta^m_t\geq B^{m,\ta}_{N^{m,\ta}_{t-}+1}}\eta^{\ta} \\
& & \hspace{5cm} + \; {\bf 1}_{\beta^0_t\geq B^{\nta}_{N^{\nta}_{t-}+1}} \eta^{\nta}+\sum_{i \in \llbracket 1,M \rrbracket}  X^{i,\beta}_{t-}\eta^{\soc}\Big)dt \Big] \\
&=&  \E\Big[\int_0^\infty H_t\frac{1}{M}\sum_{m \in \llbracket 1,M \rrbracket}  (1-X^{m, \beta}_t)\Big(\eta^{\i}+\P[b\geq B^{m,\ta}_{1}]_{b:= \beta^m_t}\eta^{\ta} \\
& & \hspace{5cm} + \; \P[b\geq B^{\nta}_{1}]_{b:=\beta^0_t}  \eta^{\nta}+\sum_{i \in \llbracket 1,M\rrbracket}X^{i,\beta}_t\eta^{\soc}\Big)dt\Big]. 
\enqs
This expression can be rewritten  as 
\begin{align}\label{eq-dp}
\E\Big[\int_0^\infty H_tdp^\beta_t\Big] &= \E\Big[\int_0^\infty H_t G_t^\beta dt \Big],
\end{align}
where $\alpha^{\beta}_t:=\frac{\sum_{m=1}^M(1-X^{m, \beta}_{t^-})\P[b\geq B^{\ta}_{1}]_{b:= \beta^m_t}}{M(1-p^\beta_{t^-})}$, and  
\begin{align}  \label{defG} 
G_t^\beta &:= \;  (1-p^\beta_t)\big(\eta^{\i}+ \eta^{\ta} \alpha^\beta_t+ \eta^{\nta} \P[b\geq B^{\nta}_{1}]_{b:=\beta^0_t} +\eta^{\soc} p^\beta_t\big), \quad \forall t\in\R_+. 
\end{align} 
This means that $G^\beta$ is the intensity process of $p^\beta$.  

\vspace{1mm}

\noindent {\it Step 2: Lower bound for $V(\beta)$.}  
From \eqref{defVsocialno}, and using the intensities of the Poisson processes,   we rewrite the cost function as 
\beqs 
V(\beta)
&=&\E\Big[\sum_{m=1}^M\int_0^\infty \Big(K(1-X^{m,\beta}_{t^-})+{\bf 1}_{\beta^m_t\geq B^{m,\ta}_{N_t^{m,\ta}}} \bc^{\ta}(\beta^m_t,B^{m,\ta}_{N_t^{m,\ta}}) \eta^{\ta} \\
& & \hspace{5cm} + \;  {\bf 1}_{\beta^0_t\geq B^{\nta}_{N^{\nta}_t}} \bc^{\nta}(\beta^0_t, B^{\nta}_{N^{\nta}_t}) \eta^{\nta}\Big)dt\Big]  \\
&=&\E\Big[\sum_{m=1}^M\int_0^\infty \Big(K(1-X^{m,\beta}_{t^-}) + \eta^{\ta} \E[\bc^{\ta}(b,B^{1,\ta}_{1}){\bf 1}_{b\geq B^{1,\ta}_{1}}]_{b:=\beta^m_t}   \\
& & \hspace{5cm} + \;  \eta^{\nta} \E[\bc^{\nta}(b, B^{\nta}_{1}){\bf 1}_{b\geq B^{\nta}_{1}}]_{b:=\beta^0_t}  \Big)dt\Big]. 
\enqs 
Now, we can bound from below the part  $\E[\bc^{\ta}(b,B^{1,\ta}_{1}){\bf 1}_{b\geq B^{1,\ta}_{1}}]_{b:=\beta^m_t}$ by $(1-X^{m,\beta}_{t^-})\E\big[\bc^{\ta}(b,B^{1,\ta}_{1}){\bf 1}_{b\geq B^{1,\ta}_{1}}\big]_{b:=\beta^m_t}$ and the part 
 $\E[\bc^{\nta}(b, B^{\nta}_{1}){\bf 1}_{b\geq B^{\nta}_{1}}]_{b:=\beta^0_t}$ by ${\bf 1}_{p^\beta_{t^-}<1}\E[\bc^{\nta}(b, B^{\nta}_{1}){\bf 1}_{b\geq B^{\nta}_{1}}]_{b:=\beta^0_t}$, so that 
\begin{align} 
V(\beta) 
&\geq \;  M\E\Big[\int_0^\infty \Big(K(1-p^{\beta}_{t^-})+\frac{1}{M}\sum_{m=1}^M(1-X^{m,\beta}_{t^-})  \eta^{\ta} \E\big[\bc^{\ta}(b,B^{1,\ta}_{1}){\bf 1}_{b\geq B^{1,\ta}_{1}}\big]_{b:=\beta^m_t} \\
&  \hspace{4cm} + \;  {\bf 1}_{p^\beta_{t^-}<1} \eta^{\nta}  \E\big[ \bc^{\nta}(b, B^{\nta}_{1}) {\bf 1}_{b\geq B^{\nta}_{1}}\big]_{b:=\beta^0_t} \Big)dt\Big] \\
&= \;  M \E\Big[\int_0^\infty H_t G_t^\beta d t\Big],   \label{inegV} 
\end{align}  
where $H_t$ $=$ $\tilde H_t(p_{t^-}^\beta)$ with  $\tilde H_t(p)$  defined for $p$ $\in$ $[0,1]$ by 
\begin{align}  
\tilde H_t(p) &:= \;   \frac{K+ \frac{\eta^{\ta}}{M(1-p)}\Sum_{m=1}^M(1-X^{m,\beta}_{t-})    \E[\bc^{\ta}(b,B^{1,\ta}_{1}){\bf 1}_{b\geq B^{1,\ta}_{1}}]_{b:=\beta^m_t} 
+ {\bf 1}_{p<1}  \eta^{\nta} \E\big[ \frac{\bc^{\nta}(b, B^{\nta}_{1})}{1-p}  {\bf 1}_{b\geq B^{\nta}_{1}}\big]_{b:=\beta^0_t} }{\eta^{\i}+ \eta^{\ta} \alpha^\beta_t+ \eta^{\nta} \P[b\geq B^{\nta}_{1}]_{b:=\beta^0_t} + \eta^{\soc} p },
\end{align} 
(with the convention that $\frac{0}{0}=0$) by recalling  the definition of $G^\beta$ in \eqref{defG}. 
For such $H$, which is clearly a positive and $\F$-predictable bounded process, we have from \eqref{inegV} and \eqref{eq-dp} 
\beqs 
V(\beta)&\geq& 
M \E\Big[\int_0^\infty \tilde H_t (p_{t^-}^\beta) dp^{\beta}_t\Big].    
\enqs 
The above r.h.s. is turned into a sum over successive values of $p^{\beta}_t$ as
\begin{align}
V(\beta) & \geq \; \E\Big[\sum_{p\in \P_M}  \tilde H_{\tau_p^\beta}(p) \Big],  
\end{align} 
where $\tau^\beta_p$ $=$ $\inf\{t\in\R_+: p^\beta_t=p+1/M\}$  for $p\in \P_M$.  
Notice that in $\tilde H_{\tau_p^\beta}(p)$,  the terms $\sum_{m=1}^M(1-X^{m,\beta}_{\tau^{\beta}_p-})\E[\bc^{\ta}(b^{m,\ta},B^{1,\ta}_{1}){\bf 1}_{b\geq B^{1,\ta}_{1}}]_{b:=\beta^m_{\tau^{\beta}_p}}$ and the sum in the definition of $\alpha^\beta_{\tau^{\beta}_p}$ are only summing over the $M(1-p)$ indices $m$ such that 
$X^{m,\beta}_{\tau^{\beta}_p-}=0$. We can thus clearly bound from below $\tilde H_{\tau_p^\beta}(p)$, for $p$ $\in$ $\P_M$, by 
\beqs
\tilde H_{\tau_p^\beta}(p) & \geq & \inf_{\underset{m\in\llbracket 1, M(1-p)\rrbracket}{b^{m,\ta}, b^{\nta}\in \R_+}}  w(p;(b^{m,\ta})_{m\in\llbracket 1,M(1-p)\rrbracket},b^{\nta}),  
\enqs
with
\beqs
& & w(p;(b^{m,\ta})_{m\in\llbracket 1,M(1-p)\rrbracket},b^{\nta}) \\
&:=& \frac{K+\frac{\eta^{\ta}}{M(1-p)}\Sum_{m=1}^{M(1-p)}\E\big[\bc^{\ta}(b^{m,\ta},B^{1,\ta}_{1}){\bf 1}_{b^{m,\ta}\geq B^{1,\ta}_{1}}\big]  
+ \eta^{\nta} \E\big[\frac{\bc^{\nta}(b^{\nta}, B^{\nta}_{1})}{1-p}{\bf 1}_{b^{\nta}\geq B^{\nta}_{1}}\big]  }{\eta^{\i}+  \eta^{\ta} \frac{\Sum_{m=1}^{M(1-p)}\P[b^{m,T}\geq B^{\ta}_{1}]}{M(1-p)} + \eta^{\nta} \P[b^{\nta}\geq B^{\nta}_{1}] +\eta^{\soc}p }, 
\enqs
so that for any open-loop bidding map control $\beta$ $\in$ $\Pi_{OL}$, 
\begin{align} \label{Vlower} 
V(\beta) & \geq \;  \sum_{p\in \P_M} v(p),
\end{align} 
with 
\begin{align}  \label{eq-infimum}
v(p) &:= \;  \inf_{\underset{m\in\llbracket 1, M(1-p)\rrbracket}{b^{m,\ta}, b^{\nta}\in \R_+}}  w(p;(b^{m,\ta})_{m\in\llbracket 1,M(1-p)\rrbracket},b^{\nta}), \quad p \in \P_M.  
\end{align} 

\vspace{1mm}

\noindent {\it Step 3:  Attaining  the lower bound}.  By definition of $v(p)$ in \eqref{eq-infimum}, we have for all $b^{m,\ta}\in\R_+$, $m\in\llbracket 1, M(1-p)\rrbracket$, and  $b^{\nta}\in \R_+$,  
\beqs 
&&K+\frac{ \eta^{\ta}}{M(1-p)}\sum_{m=1}^{M(1-p)}\E[\bc^{\ta}(b^{m,\ta},B^{1,\ta}_{1}){\bf 1}_{b^{m,\ta}\geq B^{1,\ta}_{1}}] + \eta^{\nta}  \E[\frac{\bc^{\nta}(b^{\nta}, B^{\nta}_{1})}{1-p}{\bf 1}_{b^{\nta}\geq B^{\nta}_{1}}]  \\
&\geq& v(p)\Big(\eta^{\i}+ \eta^{\ta} \frac{\sum_{m=1}^{M(1-p)}\P[b^{m,\ta}\geq B^{\ta}_{1}]}{M(1-p)}+ \eta^{\nta}\P[b^{\nta}\geq B^{\nta}_{1}]  +  \eta^{\soc} p \Big)
\enqs 
which is  equivalent to
\beqs 
&&K-v(p) (\eta^{\i}+p\eta^{\soc})+\frac{\eta^{\ta} }{M(1-p)}\sum_{m=1}^{M(1-p)}\E\big[\big(\bc^{\ta}(b^{m,\ta},B^{1,\ta}_{1})-v(p)\big){\bf 1}_{b^{m,\ta}\geq B^{1,\ta}_{1}}\big]   \\
&& \quad + \;  \eta^{\nta} \E\Big[\Big(\frac{\bc^{\nta}(b^{\nta}, B^{\nta}_{1})}{1-p}-v(p)\Big){\bf 1}_{b^{\nta} \geq B^{\nta}_{1}}\Big]   \; \geq \; 0. 
\enqs 
Moreover, we have equality if and only if $b^{m,\ta}\in\R_+$, $m\in\llbracket 1, M(1-p)\rrbracket$, and  $b^{\nta}\in \R_+$ attains the infimum in \eqref{eq-infimum}. This means  that
\begin{align}
& \quad \;\argmin_{\underset{m\in\llbracket 1, M(1-p)\rrbracket}{b^{m,\ta}, b^{\nta}\in \R_+}}  w(p;(b^{m,\ta})_{m\in\llbracket 1,M\rrbracket},b^{\nta}) \\
&\;=\argmin_{\underset{m\in\llbracket 1, M(1-p)\rrbracket}{b^{m,\ta}, b^{\nta}\in \R_+}} \Big\{  
\frac{\eta^{\ta}}{M(1-p)}\sum_{m=1}^{M(1-p)}\E[(\bc^{\ta}(b^{m,\ta},B^{1,\ta}_{1})-v(p)){\bf 1}_{b^{m,\ta}\geq B^{1,\ta}_{1}}]  \\
& \hspace{5cm}  + \;   \eta^{\nta} \E\Big[\Big(\frac{\bc^{\nta}(b^{\nta} , B^{\nta}_{1})}{1-p}-v(p)\Big){\bf 1}_{b^{\nta} \geq B^{\nta}_{1}}\Big] \Big\} \\
&\;=\Big(\Prod_{m\in\llbracket 1, M(1-p)\rrbracket}\argmin_{b^{m,\ta}\in \R_+}\E[(\bc^{\ta}(b^{m,\ta},B^{1,\ta}_{1})-v(p)){\bf 1}_{b^{m,\ta}\geq B^{1,\ta}_{1}}]\Big) \\
& \hspace{3cm}  \times \argmin_{b^{\nta}\in \R_+}  \E\Big[\Big(\frac{\bc^{\nta}(b^{\nta} , B^{\nta}_{1})}{1-p}-v(p)\Big){\bf 1}_{b^{\nta} \geq B^{\nta}_{1}}\Big] \\
&\;= \; \Big(\argmax_{b^{\ta}\in \R_+}\E\big[ \big( v(p) - \bc^{\ta}(b^{\ta},B^{1,\ta}_{1}) \big) {\bf 1}_{b^{\ta}\geq B^{1,\ta}_{1}}\big]\Big)^{M(1-p)} \\
& \hspace{3cm}  \times \argmax_{b^{\nta}\in \R_+}\E\Big[\Big(v(p)-\frac{\bc^{\nta}(b^{\nta} , B^{\nta}_{1})}{1-p}\Big){\bf 1}_{b^{\nta} \geq B^{\nta}_{1}}\Big].   \label{eq-split-static}
\end{align}
Relation \eqref{eq-split-static} shows that the original argmin over $M(1-p)+1$ arguments in $v(p)$  is reduced into the search of two  argmax for the optimal bid of static auctions, one with respect to  maximal bid  distribution $\Lc(B^{1,\ta}_{1})$ for the targeted auction, and the other  with respect to the 
maximal bid distribution $\Lc(B^{\nta}_{1})$ of the non-targeted auction. 
 
Let us check that these sets are not empty. We study the set
\begin{align} \label{setmin}  
\argmax_{b^{\ta}\in \R_+}\E[(v(p)-\bc^{\ta}(b^{\ta},B^{1,\ta}_{1})){\bf 1}_{b^{\ta}\geq B^{1,\ta}_{1}}],
\end{align} 
the other set being treated similarly, and distinguish the two cases of paying auction rules: 
\begin{itemize}
\item[1.]  {\it First-price auction: $\bc^{\ta}(b,B)=b$}. In this case, the set in \eqref{setmin} is written as 
\beqs 
\argmax_{b^{\ta}\in \R_+}  \big\{ (v(p)-b^{\ta}) \P[ b^{\ta}\geq B^{1,\ta}_{1} ] \big\}. 
\enqs 
Notice that for $b^{\ta}>v(p)$, the expression $(v(p)-b^{\ta})\P[ b^{\ta}\geq B^{1,\ta}_{1} ]$ is strictly negative, and thus smaller than its value for  $b^{\ta}=0$. The maximisation can thus be restricted to $[0,v(p)]$. Notice that $b^{\ta}\mapsto v(p)-b^{\ta}$ is positive and continuous on $[0,v(p)]$, and thus upper semi-continuous, and $b^{\ta}\mapsto \P[b^{\ta}\geq B^{1,\ta}_{1}]$ is positive, non-decreasing and càd-làg, and thus upper semi-continuous. 
It follows that $b^{\ta}\mapsto (v(p)-b^{\ta})\P[b^{\ta}\geq B^{1,\ta}_{1}]$ is positive and upper semi-continuous on $[0,v(p)]$, and thus reaches its maximum.
\item[2.] {\it Second-price auction case:  $\bc^{\ta}(b,B)=B$}.  In this case, the set  \eqref{setmin} is written as 
\beqs 
\argmax_{b^{\ta}\in \R_+}\E[(v(p)-B^{1,\ta}_{1}){\bf 1}_{b^{\ta}\geq B^{1,\ta}_{1}}],
\enqs 
and it is clear that the maximum is attained  by $b^{\ta}=v(p)$. 
\end{itemize}
This proves that the set \eqref{eq-split-static} is not empty. Moreover, any $b^{m,\ta}\in\R_+$, $m\in\llbracket 1, M(1-p)\rrbracket$, and  $b^{\nta}\in \R_+$ in the set \eqref{eq-split-static} reaches the infimum in \eqref{eq-infimum}. Given the form of the set \eqref{eq-split-static}, one can clearly take an element of the form 
$((b^{\ta})_{m\in\llbracket 1, M(1-p)\rrbracket}, b^{\nta})$ in this set (i.e. the same bid $b^{\ta}$ for the targeted advertising bid associated to all individuals who are not informed yet). Thus, the infimum in \eqref{eq-infimum} is written as 
\begin{align}
& \quad v(p) \\
& = \;\inf_{b^{\ta}, b^{\nta}\in \R_+}\frac{K+\frac{ \eta^{\ta}}{M(1-p)}\sum_{i=1}^{M(1-p)}\E[\bc^{\ta}(b^{\ta},B^{1,\ta}_{1}){\bf 1}_{b^{\ta}\geq B^{1,\ta}_{1}}] + \eta^{\nta} \E[\frac{\bc^{\nta}(b^{\nta} , B^{\nta}_{1})}{1-p}{\bf 1}_{b^{\nta} \geq B^{\nta}_{1}}]  }
{\eta^{\i}+ \eta^{\ta}  \frac{\sum_{i=1}^{M(1-p)}\P[b^{\ta}\geq B^{\ta}_{1}]}{M(1-p)} + \eta^{\nta}  \P[b^{\nta}\geq B^{\nta}_{1}] + \eta^{\soc} p }  \\
&\;=\inf_{b^{\ta}, b^{\nta}\in \R_+} v^{b^{\ta},b^{\nta}}(p),   \label{vbtbnt}
\end{align}
with
\begin{align} 
v^{b^{\ta},b^{\nta}}(p) & \; := 
\frac{K+   \eta^{\ta} \E[\bc^{\ta}(b^{\ta},B^{1,\ta}_{1}){\bf 1}_{b^{\ta}\geq B^{1,\ta}_{1}}] +  \eta^{\nta} \E[\frac{\bc^{\nta}(b^{\nta} , B^{\nta}_{1})}{1-p}{\bf 1}_{b^{\nta} \geq B^{\nta}_{1}}]  }{\eta^{\i}+ \eta^{\ta} \P[b^{\ta}\geq B^{\ta}_{1}] + \eta^{\nta} \P[b^{\nta}\geq B^{\nta}_{1}] +\eta^{\soc}p }. 
\end{align}
Therefore, we have by \eqref{Vlower} 
\begin{align}
V^\star \; = \; \inf_{\beta\in\Pi_{OL}} V(\beta) & \geq \;   \sum_{p\in \P_M}  v(p) \; = \; \sum_{p\in \P_M} \inf_{b^{\ta}, b^{\nta}\in \R_+} v^{b^{\ta},b^{\nta}}(p). 
\end{align}
Now, by considering  the control $\beta^{\mathfrak{b^\star}}$ associated to the proportion-based bidding policy $\mfb^\star$ $=$ $(\mfb^{\star,\ta},\mfb^{\star,\nta})$ defined by  
\begin{align} \label{boptim} 
(\mathfrak{b}^{\ta}(p),\mathfrak{b}^{\nta}(p)) & \in \;  \argmin_{b^{\ta}, b^{\nta}\in \R_+} v^{b^{\ta},b^{\nta}}(p),
\end{align} 
and retracing the above derivations, we see that all the inequalities turn into equalities. More precisely, the first inequality \eqref{inegV} becomes an equality whenever the bidding control used comes from a proportion-based policy. Indeed, in this case,
\beqs 
\E[\bc^{\ta}(b,B^{1,\ta}_{1}){\bf 1}_{b\geq B^{1,\ta}_{1}}]_{b:=(\beta^{\mfb})^m_t}&=&\E[\bc^{\ta}(b,B^{1,\ta}_{1}){\bf 1}_{b\geq B^{1,\ta}_{1}}]_{b:=(1-X^{m,\beta^{\mfb}}_{t^-})(\beta^{\mfb})^m_t}\\
&=&(1-X^{m,\beta^{\mfb}}_{t^-})\E\big[\bc^{\ta}(b,B^{1,\ta}_{1}){\bf 1}_{b\geq B^{1,\ta}_{1}}\big]_{b:=(\beta^{\mfb})^m_t}
\enqs 
and 
 \beqs 
 \E[\bc^{\nta}(b, B^{\nta}_{1}){\bf 1}_{b\geq B^{\nta}_{1}}]_{b:=(\beta^{\mfb})^0_t}&=&\E[\bc^{\nta}(b, B^{\nta}_{1}){\bf 1}_{b\geq B^{\nta}_{1}}]_{b:={\bf 1}_{p^{\beta^{\mfb}}_{t^-}<1}(\beta^{\mfb})^0_t}\\
 &=&{\bf 1}_{p^{\beta^{\mfb}}_{t^-}<1}\E[\bc^{\nta}(b, B^{\nta}_{1}){\bf 1}_{b\geq B^{\nta}_{1}}]_{b:=(\beta^{\mfb})^0_t},
 \enqs
 The next steps of the proof thus lead to the equality
 \beqs 
 V(\beta^{\mathfrak{b^\star}}) & = \; \E\Big[\sum_{p\in \P_M}  \tilde H_{\tau_p^{\beta^{\mathfrak{b^\star}}}}(p) \Big],
 \enqs 
 and the only other lower bound performed is $\tilde H_{\tau_p^\beta}(p)\geq v(p)$, but this bound is clearly, by definition of $\mathfrak{b^\star}$, reached by $\beta^{\mathfrak{b^\star}}$, i.e. we have $\tilde H_{\tau_p^{\beta^{\mathfrak{b^\star}}}}(p)= v(p)$, and thus $V(\beta^{\mfb^\star})$  $=$  $\sum_{p\in \P_M} v(p)$, which implies that $\beta^{\mfb^\star}$ is an optimal bidding map control:  $V(\beta^{\mfb^\star})$ $=$ $V^\star$.

\vspace{3mm}

\noindent {\it Case of a second-price auction rule, i.e.  $\bc^{\ta}(b,B)=\bc^{\nta}(b,B)=B$}.  
In this case,  relation \eqref{eq-split-static} is written as 
\begin{align}
& \quad  \argmin_{\underset{m\in\llbracket 1, M(1-p)\rrbracket}{b^{m,\ta}, b^{\nta}\in \R_+}}  w(p;(b^{m,\ta})_{m\in\llbracket 1,M\rrbracket},b^{\nta}) \\
&\;= \; \Big(\argmax_{b^{\ta}\in \R_+}\E\big[ \big( v(p) -  B^{1,\ta}_{1}) \big) {\bf 1}_{b^{\ta}\geq B^{1,\ta}_{1}}\big]\Big)^{M(1-p)} \\
& \hspace{3cm}  \times \argmax_{b^{\nta}\in \R_+}\E\Big[\Big(v(p)- \frac{   B^{\nta}_{1}}{1-p}\Big){\bf 1}_{b^{\nta} \geq B^{\nta}_{1}}\Big].   \label{eq-split-static-second}
\end{align} 
We then notice that the element $((v(p))_{m\in\llbracket 1, M(1-p)\rrbracket}, (1-p)v(p))$,  
belongs to the set \eqref{eq-split-static-second}. It follows that the infimum  in \eqref{eq-infimum} (or \eqref{vbtbnt})  can be reduced to a single parameter optimisation problem, namely
\beqs
v(p) &=& \inf_{b\in\R_+} v^{b,(1-p)b}(p).  
\enqs
Moreover, we obtain an optimal bidding map control $\beta^{\mfb^\star}$ associated to a proportion-based bidding policy in the form $\mfb^\star(p)$ $=$ $(\mrb^\star(p),(1-p)\mrb^\star(p))$, $p$ $\in$ $[0,1]$ with 
\beqs
\mrb^\star(p) & \in & \argmin_{b\in \R_+} v^{b,(1-p)b}(p). 
\enqs 

\subsubsection{Proof of the well-posedness of Definition \ref{def-smallest-proportion-based-bid}}  \label{sec:proof-def-smallest-proportion-based-bid}

Let us first prove that 
\begin{align}
& \quad \argmin_{b^{\ta}, b^{\nta}\in \R_+}v^{\bb^{\ta},\bb^{\nta}}(p)\\
= & \;\;  \argmax_{b^{\ta}\in \R_+}\E\big[(v(p)-\bc^{\ta}(b^{\ta},B^{1,\ta}_{1})){\bf 1}_{b^{\ta}\geq B^{1,\ta}_{1}}\big] \times \argmax_{b^{\nta}\in \R_+}\E\Big[\Big(v(p)-\frac{\bc^{\nta}(b^{\nta} , B^{\nta}_{1})}{1-p}\Big){\bf 1}_{b^{\nta} \geq B^{\nta}_{1}}\Big]. 
\end{align}
By definition, for any $b^{\ta},b^{\nta}\in \R_+$, we have $v^{b^{\ta},b^{\nta}}(p)$ $\geq$ $v(p)$,  with equality if and only if $b^{\ta},b^{\nta}$ reach the infimum in the definition of $v(p)$. This is formulated as 
\beqs 
\frac{K+   \eta^{\ta} \E[\bc^{\ta}(b^{\ta},B^{1,\ta}_{1}){\bf 1}_{b^{\ta}\geq B^{1,\ta}_{1}}] +  \eta^{\nta} \E[\frac{\bc^{\nta}(b^{\nta} , B^{\nta}_{1})}{1-p}{\bf 1}_{b^{\nta} \geq B^{\nta}_{1}}]  }{\eta^{\i}+ \eta^{\ta} \P[b^{\ta}\geq B^{\ta}_{1}] + \eta^{\nta} \P[b^{\nta}\geq B^{\nta}_{1}] +\eta^{\soc}p } 
&\geq& v_{}(p),
\enqs 
which is written equivalently as
\begin{align} 
K-(\eta^{\i}+\eta^{\soc}p)v_{}(p)  + \eta^{\ta} \E\big[(\bc^{\ta}(b^{\ta},B^{1,\ta}_{1})-v_{}(p)){\bf 1}_{b^{\ta}\geq B^{1,\ta}_{1}}\big]  & \\
 \quad + \; \;    \eta^{\nta} \E\Big[\Big(\frac{\bc^{\nta}(b^{\nta} , B^{\nta}_{1})}{1-p}-v_{}(p)\Big){\bf 1}_{b^{\nta} \geq B^{\nta}_{1}}\Big] &  \geq \;  0,  \label{eq-static-det0}
\end{align}
again with equality if and only if $b^{\ta},b^{\nta}$ reach the infimum in the definition of $v_{}(p)$. This clearly means  that $b^{\ta},b^{\nta}$ reach the infimum in the definition of $v_{}(p)$ if and only if they minimize \eqref{eq-static-det0} over $b^{\ta},b^{\nta}\in \R_+$, i.e. if and only if 
\begin{align}
(b^{\ta},b^{\nta}) &\in \; \argmax_{b^{\ta}\in \R_+}\E[(v(p)-\bc^{\ta}(b^{\ta},B^{1,\ta}_{1})){\bf 1}_{b^{\ta}\geq B^{1,\ta}_{1}}]  \\
& \quad \times \argmax_{b^{\nta}\in \R_+}\E\Big[\Big(v(p)-\frac{\bc^{\nta}(b^{\nta} , B^{\nta}_{1})}{1-p}\Big){\bf 1}_{b^{\nta} \geq B^{\nta}_{1}}\Big].   \label{eq-split-det}
\end{align}
It is thus clear, that there exists a unique proportion-based policy $\mfb_{min}^{\star}$, defined by
\begin{align} 
\mfb_{min}^{\star, \ta}(p)&= \; \min \argmax_{b^{\ta}\in \R_+}\E[(v(p)-\bc^{\ta}(b^{\ta},B^{1,\ta}_{1})){\bf 1}_{b^{\ta}\geq B^{1,\ta}_{1}}]  \label{bstarTmin} \\
\mfb_{min}^{\star, \nta}(p)&= \; \min \argmax_{b^{\nta}\in \R_+}\E\Big[\Big(v(p)-\frac{\bc^{\nta}(b^{\nta} , B^{\nta}_{1})}{1-p}\Big){\bf 1}_{b^{\nta} \geq B^{\nta}_{1}}\Big], \label{bstarNTmin} 
\end{align}  
such that $(\mfb_{min}^{\star, \ta}(p),\mfb_{min}^{\star, \nta}(p))$ reaches the infimum in the definition of $v(p)$, and that for any other $(\bb^{\star, \ta},\bb^{\star, \nta})$ reaching this infimum, we have
$\mfb_{min}^{\star, \ta}(p)\leq \bb^{\star, \ta}(p)$,  $\mfb_{min}^{\star, \nta}(p)\leq \bb^{\star, \nta}(p)$.

\subsubsection{Proof of Proposition \ref{propsocialno-sensitivity}}  \label{sec:prosen}

It is clear from the formula that for all $p\in \P_M$, $v(p)$ is decreasing in $\eta^{\i}$ and $\eta^{\soc}$. Regarding the sensitivity in $\eta^{\ta}$ and $\eta^{\nta}$, let us first prove that we have
\begin{align} \label{tildevp} 
v(p) &= \; \inf_{\bb^{\ta},\bb^{\nta}\in L(\Omega, \R_+), \independent B^{1,\ta}_1,B^{1,\nta}_1} v^{\bb^{\ta},\bb^{\nta}}(p),
\end{align} 
i.e. that the $\inf$ can be taken over  the set of  random variables $\bb^{\ta},\bb^{\nta}\in L(\Omega, \R_+), \independent B^{1,\ta}_1,B^{1,\nta}_1$ instead of over  $b^{\ta},b^{\nta}\in \R_+$, without changing the infimum. Denote by $\tilde v(p)$ the right hand side of \eqref{tildevp}, and 
let us check that $\tilde{v}(p)$ $=$ $v(p)$. By definition, for any $\bb^{\ta},\bb^{\nta}\in L(\Omega, \R_+), \independent B^{1,\ta}_1,B^{1,\nta}_1$, 
we have $v^{\bb^{\ta},\bb^{\nta}}(p)$ $\geq$ $\tilde v(p)$,  with equality if and only if $\bb^{\ta},\bb^{\nta}$ reach the infimum in the definition of $\tilde v(p)$. This is formulated as 
\begin{align} 
K-(\eta^{\i}+\eta^{\soc}p)v_{}(p)  + \eta^{\ta} \E\big[(\bc^{\ta}(\bb^{\ta},B^{1,\ta}_{1})-v_{}(p)){\bf 1}_{\bb^{\ta}\geq B^{1,\ta}_{1}}\big]  & \\
 \quad + \; \;    \eta^{\nta} \E\Big[\Big(\frac{\bc^{\nta}(\bb^{\nta} , B^{\nta}_{1})}{1-p}-v_{}(p)\Big){\bf 1}_{\bb^{\nta} \geq B^{\nta}_{1}}\Big] &  \geq \;  0,  \label{eq-static-rand}
\end{align}
with equality if and only if $\bb^{\ta},\bb^{\nta}$ reach the infimum in the definition of $\tilde v_{}(p)$. This clearly means  that $\bb^{\ta},\bb^{\nta}$ reach the infimum in the definition of $\tilde v_{}(p)$ if and only if they minimize \eqref{eq-static-rand} over  
$\bb^{\ta},\bb^{\nta}\in L(\Omega, \R_+), \independent B^{1,\ta}_1,B^{1,\nta}_1$. By conditioning, it is clear that \eqref{eq-static-rand} will reach the same infimum if minimized over $b^{\ta},b^{\nta}\in\R_+$. This in particular implies that the infimum in $\tilde{v}(p)$ will be reached if it is only taken over $b^{\ta},b^{\nta}\in\R_+$, 
which means that  $\tilde{v}(p)=v(p)$.

In the sequel, we stress the dependence of $v^{\bb^{\ta},\bb^{\nta}}(p)$ and $v(p)$ in $\eta^{\ta}$ by writing $v_{\eta^{\ta}}^{\bb^{\ta},\bb^{\nta}}(p)$ and $v_{\eta^{\ta}}(p)$. Let us consider $\tilde{\eta}^{\ta}<\eta^{\ta}$ and let us denote by $Z$ a Bernoulli random variable with parameter $\frac{\tilde{\eta}^{\ta}}{\eta^{\ta}}$, independent from $(B^{1,\ta}_1,B^{1,\nta}_1)$. For any $b^{\ta},b^{\nta}\in\R_+$, we define $\bb^{\ta}=Zb^{\ta}$ and $\bb^{\nta}=b^{\nta}$. Notice then that,
\beqs 
v_{\eta^{\ta}}^{\bb^{\ta},\bb^{\nta}}(p)=v_{\tilde{\eta}^{\ta}}^{b^{\ta},b^{\nta}}(p)
\enqs 
Then, we have
\beqs 
v_{\eta^{\ta}}(p)&=&\inf_{\bb^{\ta},\bb^{\nta}\in L(\Omega, \R_+), \independent B^{1,\ta}_1,B^{1,\nta}_1} v_{\eta^{\ta}}^{\bb^{\ta},\bb^{\nta}}(p)\leq \inf_{b^{\ta},b^{\nta}\in \R_+} v_{\eta^{\ta}}^{Zb^{\ta},b^{\nta}}(p)\\
&=& \inf_{b^{\ta},b^{\nta}\in \R_+} v_{\tilde{\eta}^{\ta}}^{b^{\ta},b^{\nta}}(p)= v_{\tilde{\eta}^{\ta}}(p),
\enqs 
and thus $v_{\eta^{\ta}}(p)\leq v_{\tilde{\eta}^{\ta}}(p)$, which means that $v_{\eta^{\ta}}(p)$ is decreasing in $\eta^{\ta}$. The same argument allows to prove that $v(p)$ is decreasing in $\eta^{\nta}$. 
To summarize, $v(p)$ (and thus $\inf_{\beta\in \Pi_{OL}}V(\beta)$) is decreasing in all the model's intensity parameters $\theta=(\eta^{\i},\eta^{\soc},\eta^{\ta},\eta^{\nta})$. 



Let us now study the monotonicity of the smallest optimal bid with respect to $\theta$.  We  then stress the dependence of $v^{b^{\ta},b^{\nta}}(p)$, $v(p)$, $\mfb_{min}^{\star,\ta}(p)$ and $\mfb_{min}^{\star,\nta}(p)$ in $\theta=(\eta^{\i},\eta^{\soc},\eta^{\ta},\eta^{\nta})$ by writing $v_{\theta}^{b^{\ta},b^{\nta}}(p)$, $v_{\theta}(p)$, $\mfb_{min,\theta}^{\star,\ta}(p)$ and $\mfb_{min,\theta}^{\star,\nta}(p)$. 
Let us now consider $\tilde{\theta}=(\tilde{\eta}^{\i},\tilde{\eta}^{\soc},\tilde{\eta}^{\ta},\tilde{\eta}^{\nta})$ such that $\tilde{\eta}^{\i}\leq \eta^{\i}$, $\tilde{\eta}^{\soc}\leq \eta^{\soc}$, $\tilde{\eta}^{\ta}\leq \eta^{\ta}$, $\tilde{\eta}^{\nta}\leq \eta^{\nta}$. We then know from above that  $v_{\theta}(p)\leq v_{\tilde{\theta}}(p)$. 
Let us then prove that  $\mfb_{min,\theta}^{\star,\ta}(p)\geq \mfb_{min,\tilde{\theta}}^{\star,\ta}(p)$. Assume on the contrary  that $\mfb_{min,\theta}^{\star,\ta}(p)< \mfb_{min,\tilde{\theta}}^{\star,\ta}(p)$. 
From \eqref{bstarTmin},  we have in particular that 
$\mfb_{min,\theta}^{\star,\ta}(p)\in \argmax_{b^{\ta}\in \R_+}\E[(v_\theta(p)-\bc^{\ta}(b^{\ta},B^{1,\ta}_{1})){\bf 1}_{b^{\ta}\geq B^{1,\ta}_{1}}]$, 
and thus for $b^{\ta}=\mfb_{min,\tilde{\theta}}^{\star,\ta}(p)$, we obtain
\begin{align}\label{ineq-1}
& \E[(v_{\theta}(p)-\bc^{\ta}(\mfb_{min,\theta}^{\star,\ta}(p),B^{1,\ta}_{1})){\bf 1}_{\mfb_{min,\theta}^{\star,\ta}(p)\geq  B^{1,\ta}_{1}} \\
\geq & \; \E[(v_{\theta}(p)-\bc^{\ta}(\mfb_{min,\tilde{\theta}}^{\star,\ta}(p),B^{1,\ta}_{1})){\bf 1}_{\mfb_{min,\tilde{\theta}}^{\star,\ta}(p)\geq B^{1,\ta}_{1}}].
\end{align}
On the other hand,  since  
$\mfb_{min,\tilde{\theta}}^{\star,\ta}(p)=\min \argmax_{b^{\ta}\in \R_+}\E[(v_{\tilde{\theta}}(p)-\bc^{\ta}(b^{\ta},B^{1,\ta}_{1})){\bf 1}_{b^{\ta}\geq B^{1,\ta}_{1}}]$,  by \eqref{bstarTmin}, 
we have that for all $b'^{\ta}<\mfb_{min,\tilde{\theta}}^{\star,\ta}(p)$,
\beqs 
b'^{\ta}\not\in \argmax_{b^{\ta}\in \R_+}\E[(v_{\tilde{\theta}}(p)-\bc^{\ta}(b^{\ta},B^{1,\ta}_{1})){\bf 1}_{b^{\ta}\geq B^{1,\ta}_{1}}]. 
\enqs 
Therefore, by taking  $b'^{\ta}=\mfb_{min,\theta}^{\star,\ta}(p)$ $<$  $\mfb_{min,\tilde{\theta}}^{\star,\ta}(p)$, we get 
\begin{align} \label{ineq-2}
& \E[(v_{\tilde{\theta}}(p)-\bc^{\ta}(\mfb_{min,\theta}^{\star,\ta}(p),B^{1,\ta}_{1})){\bf 1}_{\mfb_{min,\theta}^{\star,\ta}(p)\geq B^{1,\ta}_{1}}] \\
< \; & \;   \E[(v_{\tilde{\theta}}(p)-\bc^{\ta}(\mfb_{min,\tilde{\theta}}^{\star,\ta}(p),B^{1,\ta}_{1})){\bf 1}_{\mfb_{min,\tilde{\theta}}^{\star,\ta}(p)\geq B^{1,\ta}_{1}}]. 
\end{align}
By subtracting \eqref{ineq-1} to \eqref{ineq-2}, we obtain 
\beqs 
\E[(v_{\tilde{\theta}}(p)-v_{\theta}(p)){\bf 1}_{\mfb_{min,\theta}^{\star,\ta}(p)\geq B^{1,\ta}_{1}}] &<&  \E[(v_{\tilde{\theta}}(p)-v_{\theta}(p)){\bf 1}_{\mfb_{min,\tilde{\theta}}^{\star,\ta}(p)\geq B^{1,\ta}_{1}}],
\enqs 
and thus $\P[\mfb_{min,\theta}^{\star,\ta}(p)\geq B^{1,\ta}_{1}] > \P[\mfb_{min,\tilde{\theta}}^{\star,\ta}(p)\geq B^{1,\ta}_{1}]$, which contradicts the inequality $\mfb_{min,\theta}^{\star,\ta}(p)< \mfb_{min,\tilde{\theta}}^{\star,\ta}(p)$. 
This shows that  $\mfb_{min,\theta}^{\star,\ta}(p)\geq \mfb_{min,\tilde{\theta}}^{\star,\ta}(p)$. The same arguments applies to prove that $\mfb_{min,\theta}^{\star,\nta}(p)\geq \mfb_{min,\tilde{\theta}}^{\star,\nta}(p)$. 

\vspace{1mm}

Let us now study the variations of the smallest optimal bid w.r.t.  the proportion of informed people $p$. By definition of $v(p)$, we have the following properties:
\begin{itemize}
    \item If there is only targeted advertising ($\eta^{\nta}=0$), then $v(p)$ is decreasing in $p$, and the above argument applied to two proportions $\tilde{p}<p$ (instead of two model parameters $\theta,\tilde{\theta}$)  show that 
    $\bb^{\star,\ta}_{\min}(\tilde{p})\geq \mfb_{min}^{\star,\ta}(p)$.
    \item If there is no social interactions ($\eta^{\soc}=0$), $v(p)$ is increasing in $p$, and the above argument applied to two proportions $\tilde{p}<p$ (instead of two model parameters $\theta,\tilde{\theta}$) show that 
    $\bb^{\star,\ta}_{\min}(\tilde{p})\leq \mfb_{min}^{\star,\ta}(p)$.
\end{itemize}
By  \eqref{bstarNTmin}, we have
\beqs 
\mfb_{min}^{\star,\nta}(p)&=&\min\argmax_{b^{\nta} \in\R_+}\E\Big[\Big(v(p)-\frac{\bc^{\nta}(b^{\nta} , B^{\nta}_{1})}{1-p}\Big){\bf 1}_{b^{\nta} \geq B^{\nta}_{1}}\Big]\\
&=&\min\argmax_{b^{\nta} \in\R_+}\E\Big[\Big((1-p)v(p)-\bc^{\nta}(b^{\nta} , B^{\nta}_{1})\Big){\bf 1}_{b^{\nta} \geq B^{\nta}_{1}}\Big].
\enqs 
Notice that from the definition of $v(p)$ we see that $(1-p)v(p)$ is always decreasing in $p$. By the same argument as before applied to two proportions $\tilde{p}<p$ (instead of two model parameters $\theta,\tilde{\theta}$), and with $(1-p)v(p)$ playing the role of $v(p)$, we deduce  that $\mfb_{min}^{\star,\nta}(\tilde{p})\geq \mfb_{min}^{\star,\nta}(p)$.

\vspace{1mm}

Let us prove that  $\mfb_{min}^{\star,\ta}(p)\leq v(p)$ and $\bb^{\star,\nta}_{\min}(p)\leq v(p)$. We check this result for $\mfb_{min}^{\star,\nta}(p)$,  the other  being similarly proved. If the targeted ads are sold with first price auctions, i.e. if $\bc^{\ta}(b,B)=b$, then notice that for all $b^{\ta}>v(p)$, we have
\beqs 
\E[(v(p)-b^{\ta}){\bf 1}_{b^{\ta}\geq B^{1,\ta}_{1}}] \; < \;  0 \; \leq \; \E[(v(p)-0){\bf 1}_{0\geq B^{1,\ta}_{1}}].
\enqs 
This implies that any bid  $b^{\ta}>v(p)$ cannot be optimal, and thus that the smallest optimal bid $\mfb_{min}^{\star,\nta}(p)$ is smaller than $v(p)$. If the targeted ads are sold with second price auctions, we clearly have
\beqs 
\E[(v(p)-B^{1,\ta}_{1}){\bf 1}_{v(p)\geq B^{1,\ta}_{1}}] &=& \argmax_{b^{\ta}\in\R_+}\E[(v(p)-b^{\ta}){\bf 1}_{b^{\ta}\geq B^{1,\ta}_{1}}],
\enqs 
and thus $v(p)$ is an optimal bid, which implies that the smallest optimal bid $\mfb_{min}^{\star,\nta}(p)$ is, again, smaller than $v(p)$. In particular, given that
\beqs 
v(p)\leq v^{0,0}(p)=\frac{K}{\eta^{\i}+p\eta^{\soc}},
\enqs 
the smallest optimal bid $\mfb_{min}^{\star,\nta}(p)$ is bounded from above by $\frac{K}{\eta^{\i}+p\eta^{\soc}}$.

\subsection{Proof of results in Section \ref{secsocialdiscount}} 

\subsubsection{Proof of Theorem \ref{theo-socialdis} 
}

Let us fix an open-loop bidding control $\beta$. From \eqref{defVsocialdis}, we have 
\begin{align} 
V(\beta)&= \; \E\Big[\int_0^\infty   e^{-\rho t}(K(1-X^{\beta}_{t^-})dN^{\danger}_t+{\bf 1}_{\beta_t\geq B^{\ta}_{N_t^{\ta}}} B^{\ta}_{N_t^{\ta}} d N_t^{\ta})\Big]\\
&\geq \;  \E\Big[\int_0^\infty   e^{-\rho t} \Big(K(1-X^\beta_{t-})+(1-X^\beta_{t-}){\bf 1}_{\beta_t\geq B^{\ta}_{N^{\ta}_t}} B^{\ta}_{N^{\ta}_t} \eta^{\ta}\Big)dt\Big]\\
& = \;  \E\Big[\int_0^\infty \P[\tau> t] (1-X^\beta_{t-})\big(K+{\bf 1}_{\beta_t\geq B^{\ta}_{N^{\ta}_t}} B^{\ta}_{N^{\ta}_t} \eta^{\ta}\big)dt\Big]\\
&=\; \E\Big[\int_0^\infty  {\bf 1}_{\tau> t} (1-X^\beta_{t-})\big(K+{\bf 1}_{\beta_t\geq B^{\ta}_{N^{\ta}_t}} B^{\ta}_{N^{\ta}_t} \eta^{\ta}\big)dt\Big]\\
&= \; \E\Big[\int_0^{\tau} (1-X^\beta_{t-})\big(K+{\bf 1}_{\beta_t\geq B^{\ta}_{N^{\ta}_t}} B^{\ta}_{N^{\ta}_t} \eta^{\ta}\big)dt\Big],
\end{align} 
where we introduced an independent random time $\tau$ with exponential distribution of parameter $\rho$. 
Notice also that the  first inequality becomes an equality if the bidding control $\beta$ makes null bids once the individual is informed. 
Let us next  consider  a Poisson process $N$ with intensity $\rho$, whose  first time of jump is given by $\tau$, and independent of the other random variables, 
and denote by $\tilde X^\beta$ the process satisfying the dynamic
\beqs 
\tilde{X}^{\beta}_{0^-}&=& 0\\
d\tilde{X}^{\beta}_t&=&(1-\tilde{X}^\beta_{t^-})(dN^{\i}_t+dN_t+{\bf 1}_{\beta_t \geq B^{\ta}_{N_t^{\ta}}}dN^{\ta}_t), \quad t \geq 0. 
\enqs
Notice that $\tilde{X}^{\beta}$ has exactly the same dynamic as $X^\beta$ except that there is an additional cause of transition to state $1$ given by the term $dN$. It is then clear that we have
\beqs 
\E\Big[\int_0^{\tau} (1-X^\beta_{t-})(K+{\bf 1}_{\beta_t\geq B^{\ta}_{N^{\ta}_t}} B^{\ta}_{N^{\ta}_t} \eta^{\ta})dt\Big] \; = \; \E\Big[\int_0^\infty (1-\tilde{X}^\beta_{t^-})(K+{\bf 1}_{\beta_t\geq B^{\ta}_{N^{\ta}_t}} B^{\ta}_{N^{\ta}_t} \eta^{\ta})dt\Big]\\
=  \;  \E\Big[\int_0^\infty (K(1-\tilde{X}^\beta_{t^-})dN^{\danger}_t+{\bf 1}_{\tilde{\beta}_t\geq B^{\ta}_{N^{\ta}_t}} B^{\ta}_{N^{\ta}_t} dN^{\ta}_t)\Big], 
\enqs 
where $\tilde{\beta}_t=(1-\tilde{X}^\beta_{t^-})\beta_t$. By noting $\tilde{N}^{\i}=N^{\i}+N$, we obtain a Poisson process $\tilde{N}^{\i}$ with intensity $\eta^\i+\rho$, and the dynamic of $\tilde X^\beta$ is rewritten as
\beqs 
\tilde{X}^{\beta}_{0^-} &=& 0\\
d\tilde{X}^{\beta}_t&=&(1-\tilde{X}^\beta_{t^-})(d\tilde{N}^{\i}_t+{\bf 1}_{\beta_t \geq B^{\ta}_{N^{\ta}}}dN^{\ta}_t), \quad t \geq 0.
\enqs
The cost $V(\beta)$ is thus bounded from below by the cost associated to the bidding map control $\beta$ $=$ $(0,\tilde{\beta})$ of problem in Section \ref{sec:socialno}, 
with a population of  $M=1$ individual,  with an intensity $\eta^\i+\rho$ for the counting number of connections on the website with information $\i$, 
and where $\eta^{\nta}=\eta^{\soc}=0$, i.e. the individual never connects to a website displaying non-targeted ads, and individuals do not socially interact. 
From the result proved in Section \ref{sec:socialno}, we then know that  $V(\beta)$ is thus bounded from below by: 
\beqs 
V(\beta)&\geq&  \inf_{b\in\R_+}\frac{K+  \eta^{\ta}\E[B^{\ta}_1{\bf 1}_{b\geq B^{1,\ta}_1}]}{\eta^{\i}+\rho +\eta^{\ta}
\P[b\geq B^{\ta}_1]}. 
\enqs 
It is then direct  to retrace this derivation with the particular bidding control $\beta^{b_\star}$ associated to the constant bidding policy $b^\star$ such that
\beqs 
b^\star&=& \argmin_{b\in \R_+}\frac{K+  \eta^{\ta}\E[B^{\ta}_1{\bf 1}_{b\geq B^{\ta}_1}]}{\eta^{\i}+\rho +\eta^{\ta}\P[b\geq B^{\ta}_1]},
\enqs 
and to turn inequalities into equalities. This concludes the result for this case.

\subsubsection{Proof of Proposition \ref{propsocial-sensitivity}} 

Notice that the optimal value $V^\star$ corresponds to the optimal value of the problem solved in Theorem \ref{theo-social-no-discount} with $M=1$, $\eta^{\nta}=0$ and $\eta^{\i}$ replaced by $\eta^{\i}+\rho$. The sensitivity of the optimal value and smallest optimal bids to the model's parameters directly follows.

\subsection{Proof of results in Section \ref{sec:compur}} 

\subsubsection{Proof of Theorem \ref{theo-purchase} 
}
The idea is again  to reduce to the previous case. Given an open-loop bidding control $\beta$, we have by \eqref{defVpurchase} and denoting by  $\tau^\beta$ $=$ $\inf\{ t\geq 0: X_{t}^\beta =1\}$:  
\beqs 
V(\beta)&=& \E\Big[ e^{-\rho \tau^\beta}K-\int_0^\infty   e^{-\rho t}{\bf 1}_{\beta_t\geq B^{\ta}_{N^{\ta}_t}} B^{\ta}_{N^{\ta}_t} d N^{\ta}_t)\Big]\\
&=& \E\Big[\int_{\tau^\beta}^\infty \rho e^{-\rho t}K-\int_0^\infty   e^{-\rho t}{\bf 1}_{\beta_t\geq B^{\ta}_{N^{\ta}_t}} B^{\ta}_{N^{\ta}_t} d N^{\ta}_t)\Big]\\
&=& \E\Big[\int_0^\infty  e^{-\rho t}\rho KX^\beta_{t-}dt-\int_0^\infty  e^{-\rho t}{\bf 1}_{\beta_t\geq B^{\ta}_{N^{\ta}_t}} B^{\ta}_{N^{\ta}_t} d N^{\ta}_t)\Big]. 
\enqs 
The problem is thus reduced to a continuous gain problem, with continuous reward $\rho K$ from the time of information. This continuous gain problem is then turned into a continuous cost problem as follows: 
\beqs 
V(\beta) 
&=& K-\E\Big[\int_0^\infty  e^{-\rho t} (\rho K(1-X^\beta_{t-})dt+{\bf 1}_{\beta_t\geq B^{\ta}_{N^{\ta}_t}} B^{\ta}_{N^{\ta}_t}d N^{\ta}_t)\Big]\\
&=& K-\E\Big[\int_0^\infty   e^{-\rho t} (\rho K(1-X^\beta_{t-})dN^{\danger}_t+{\bf 1}_{\beta_t\geq B^{\ta}_{N^{\ta}_t}} B^{\ta}_{N^{\ta}_t}d N^{\ta}_t)\Big]. 
\enqs 
We are reduced to the previous case (social marketing with discount factor). This concludes the proof.

\subsubsection{Proof of Proposition \ref{prop-sensitivity}}

The cost dual viewpoint reduces the problem to the problem of social marketing with discount factor $\rho$ and continuous cost $\rho K$. This directly yields the sensitivity of the optimal value and smallest optimal bid in all the model's parameters except $\rho$, since here $\rho$ also appears in the continuous cost $\rho K$. For the sensitivity in $\rho$, it is suitable  to use the gain viewpoint of the optimal value, in Theorem \ref{theo-purchase}. In this expression, it is clear that $V^\star$ is decreasing in $\rho$. From this property, a similar argument as used several times in the proof of Proposition \ref{propsocialno-sensitivity} allows to conclude that the smallest optimal bid is then increasing in $\rho$. 

\subsection{Proof of Theorem \ref{theo-sus}  
}
Given an open-loop bidding control $\beta$, we have by \eqref{defVsus} 
\beqs 
V(\beta)&=&  \E\Big[\sum_{n\in \N}  e^{-\rho(\tau^\beta+n)} K - \int_0^\infty   e^{-\rho t} {\bf 1}_{\beta_t\geq B^{\ta}_{N^{\ta}_t}} B^{\ta}_{N^{\ta}_t} d N^{\ta}_t)\Big]\\
&=&  \E\Big[ e^{-\rho \tau^\beta} \frac{K}{1-e^{-\rho}}- \int_0^\infty  e^{-\rho t} {\bf 1}_{\beta_t\geq B^{\ta}_{N^{\ta}_t}} B^{\ta}_{N^{\ta}_t} d N^{\ta}_t)\Big]. 
\enqs 
This reduces the problem to the purchase-based case  studied in the previous section,  and concludes the proof.

 \section{Conclusion} \label{sec:conclusion}

In this paper, we have developed several targeted advertising models with semi-explicit solutions. An important feature  of these models is a very concrete description of the ``modern" advertising problem. 
One or several individuals are really modelled through their behaviours that  involve connections to various types of websites at random times as well as social interactions. The advertising auctions are also precisely defined by  considering various auction rules (second-price auctions, first-price auctions). Several variants of our models, that we did not study in this work for the sake of conciseness, can be easily analysed with the techniques developed in this work. For instance, in the first three models with a single Individual, one could study a model where first-price auctions and second-price auctions coexist, which would lead to more general formulas. There is also room for exploration to enrich the models while keeping them tractable with semi-explicit solutions: in the fourth model with an interacting population, it might be possible to add a bit of heterogeneity in the population connections and social interactions. It would be also interesting and useful in practice to consider the alternative  for individuals  
not to be receptive with some probability to the information (hence not purchasing a product, or continuing to behave ``dangerously"). 
Another opportune development, regarding the auctions, could be to model the maximal bid from other bidders more realistically than with an i.i.d. sequence of random variables, for instance as a Markov chain.  Another approach could be to explicitly model several bidding agents, for instance  playing according to the so-called fictitious play principle. In such game, several bidding Agents have pieces of Information to diffuse to Individuals, and each time when they declare  their bid, they follow the strategies  according to the model studied in this paper by modelling the other bidders' maximal bid as a sequence of i.i.d. random variables distributed as the empirical distribution of past maximal bids. Notice that this would require  for  Agents to constantly recalibrate their model, as new auctions modify the empirical distribution of past maximal bids.

\bibliographystyle{plain}
\bibliography{bibliography}
 
\end{document}